\documentclass[10pt,twoside, a4paper,reqno]{amsart}

\usepackage{amscd,amssymb,amsmath,graphicx,verbatim,mathrsfs,xcolor}
\usepackage[british]{babel}
\usepackage{wasysym}
\usepackage{hyperref}
\usepackage{enumitem}
\usepackage{microtype}
\usepackage[capitalize]{cleveref}
\usepackage{mathtools}
\usepackage{multicol}
\usepackage[margin=2.5cm]{geometry}
 
\usepackage[normalem]{ulem}
\usepackage[foot]{amsaddr}
\usepackage[only,llbracket,rrbracket]{stmaryrd}
\usepackage{manfnt}
\usepackage{fontawesome5}
\binoppenalty=\maxdimen
\relpenalty=\maxdimen

\usepackage[titletoc]{appendix} %
\AtBeginEnvironment{appendices}{\def\chaptername\appendixname}
\AtEndEnvironment{appendices}{\def\chaptername\oldchaptername}
{\endappendices}

\usepackage{tikz,mathrsfs}
\usetikzlibrary{babel,arrows,decorations.markings,arrows.meta,decorations.pathmorphing,cd,patterns,backgrounds}
\tikzset{
	symbol/.style={
		draw=none,
		every to/.append style={
			edge node={node [sloped, allow upside down, auto=false]{$#1$}}}
	}}
\tikzset{>=stealth'}
\tikzcdset{arrow style=tikz, diagrams={>=stealth'}, row sep=6em, column sep=6em}
\def\arrowLengthDisplayStyle{4ex}
\def\arrowHeightDisplayStyle{.8ex}
\def\arrowSkipDisplayStyle{.5ex}
\def\arrowLengthTextStyle{3ex}
\def\arrowHeightTextStyle{.8ex}
\def\arrowSkipTextStyle{.4ex}
\def\arrowLengthScriptStyle{2.5ex}
\def\arrowHeightScriptStyle{.6ex}
\def\arrowSkipScriptStyle{.3ex}
\def\arrowLengthScriptScriptStyle{2ex}
\def\arrowHeightScriptScriptStyle{.4ex}
\def\arrowSkipScriptScriptStyle{.2ex}

\renewcommand{\to}{\arrow{->}}

\newcommand{\embed}{\arrow{right hook->}}
\renewcommand{\mapsto}{\arrow{|->}}

\newcommand{\MakeTikzArrowWithSuperscriptSubscript}[4]
{
	\mathchoice
	{ 
		\hspace*{\arrowSkipDisplayStyle}
		\begin{tikzpicture}[baseline]
		\draw [#1] (0,\arrowHeightDisplayStyle) -- node [above] {$#2$} node [below] {$#3$} (#4 * \arrowLengthDisplayStyle, \arrowHeightDisplayStyle);
		\end{tikzpicture}
		\hspace*{\arrowSkipDisplayStyle}
	}
	{ 
		\hspace*{\arrowSkipTextStyle}
		\begin{tikzpicture}[baseline]
		\draw [#1] (0,\arrowHeightTextStyle) -- node [above] {$\scriptstyle #2$} node [below] {$\scriptstyle #3$} (#4 * \arrowLengthTextStyle, \arrowHeightTextStyle);
		\end{tikzpicture}
		\hspace*{\arrowSkipTextStyle}
	}
	{ 
		\hspace*{\arrowSkipScriptStyle}
		\begin{tikzpicture}[baseline]
		\draw [#1] (0,\arrowHeightScriptStyle) -- node [above] {$\scriptscriptstyle #2$} node [below] {$\scriptscriptstyle #3$} (#4 * \arrowLengthScriptStyle, \arrowHeightScriptStyle);
		\end{tikzpicture}
		\hspace*{\arrowSkipScriptStyle}
	}
	{ 
		\hspace*{\arrowSkipScriptScriptStyle}
		\begin{tikzpicture}[baseline]
		\draw [#1] (0,\arrowHeightScriptScriptStyle) -- node [above] {$\scriptscriptstyle #2$} node [below] {$\scriptscriptstyle #3$} (#4 * \arrowLengthScriptScriptStyle, \arrowHeightScriptScriptStyle);
		\end{tikzpicture}
		\hspace*{\arrowSkipScriptScriptStyle}
	}
}

\newcommand{\MakeTikzArrowWithCentralLabel}[3]
{
	\mathchoice
	{ 
		\hspace*{\arrowSkipDisplayStyle}
		\begin{tikzpicture}[baseline]
		\draw [#1] (0,\arrowHeightDisplayStyle) -- node [fill=white,inner sep=1pt] {$#2$} (#3 * \arrowLengthDisplayStyle, \arrowHeightDisplayStyle);
		\end{tikzpicture}
		\hspace*{\arrowSkipDisplayStyle}
	}
	{ 
		\hspace*{\arrowSkipTextStyle}
		\begin{tikzpicture}[baseline]
		\draw [#1] (0,\arrowHeightTextStyle) -- node [fill=white,inner sep=1pt] {$\scriptstyle #2$} (#3 * \arrowLengthTextStyle, \arrowHeightTextStyle);
		\end{tikzpicture}
		\hspace*{\arrowSkipTextStyle}
	}
	{ 
		\hspace*{\arrowSkipScriptStyle}
		\begin{tikzpicture}[baseline]
		\draw [#1] (0,\arrowHeightScriptStyle) -- node [fill=white,inner sep=1pt] {$\scriptscriptstyle #2$} (#3 * \arrowLengthScriptStyle, \arrowHeightScriptStyle);
		\end{tikzpicture}
		\hspace*{\arrowSkipScriptStyle}
	}
	{ 
		\hspace*{\arrowSkipScriptScriptStyle}
		\begin{tikzpicture}[baseline]
		\draw [#1] (0,\arrowHeightScriptScriptStyle) -- node [fill=white,inner sep=1pt] {$\scriptscriptstyle #2$} (#3 * \arrowLengthScriptScriptStyle, \arrowHeightScriptScriptStyle);
		\end{tikzpicture}
		\hspace*{\arrowSkipScriptScriptStyle}
	}
}

\def\arrow#1{\def\lastArrowStyle{#1}
	\futurelet\testchar\arrowMaybeStreched}
\def\arrowMaybeStreched{\ifx[\testchar \let\next\arrowStreched
	\else \let\next\arrowUnstreched \fi
	\next}

\def\arrowStreched[#1]{\def\lastArrowStrech{#1}
	\futurelet\testchar\arrowMaybeLabel}
\def\arrowUnstreched{\def\lastArrowStrech{1}
	\futurelet\testchar\arrowMaybeLabel}

\def\arrowMaybeLabel{\ifx^\testchar \let\next\arrowSuperscript
	\else \ifx_\testchar \let\next\arrowSubscript
	\else \ifx~\testchar \let\next\arrowCentralLabel
	\else \let\next\arrowNoLabel
	\fi
	\fi
	\fi
	\next}

\def\arrowSuperscript^#1{\def\lastArrowSuperscript{#1}
	\futurelet\testchar\arrowSuperMaybeSub}
\def\arrowSuperMaybeSub{\ifx_\testchar \let\next\arrowSuperscriptSubscript
	\else \let\next\arrowSuperscriptNoSubscript \fi
	\next}

\def\arrowSubscript_#1{\def\lastArrowSubscript{#1}
	\futurelet\testchar\arrowSubMaybeSuper}
\def\arrowSubMaybeSuper{\ifx^\testchar \let\next\arrowSubscriptSuperscript
	\else \let\next\arrowSubscriptNoSuperscript \fi
	\next}

\def\arrowSuperscriptSubscript_#1{\def\lastArrowSubscript{#1}
	\arrowDrawSupSub}
\def\arrowSuperscriptNoSubscript{\def\lastArrowSubscript{}
	\arrowDrawSupSub}
\def\arrowSubscriptSuperscript^#1{\def\lastArrowSuperscript{#1}
	\arrowDrawSupSub}
\def\arrowSubscriptNoSuperscript{\def\lastArrowSuperscript{}
	\arrowDrawSupSub}
\def\arrowNoLabel{\def\lastArrowSuperscript{}
	\def\lastArrowSubscript{}
	\arrowDrawSupSub}

\def\arrowCentralLabel~#1{\MakeTikzArrowWithCentralLabel{\lastArrowStyle}{#1}{\lastArrowStrech}}
\def\arrowDrawSupSub{\MakeTikzArrowWithSuperscriptSubscript{\lastArrowStyle}{\lastArrowSuperscript}{\lastArrowSubscript}{\lastArrowStrech}}

\tikzcdset{arrow style=tikz, diagrams={>=stealth'}, row sep=4em, column sep=4em}

\counterwithin*{equation}{subsection}

\makeatletter

\renewcommand\subsection{\@secnumfont}{\bfseries}%
\renewcommand\subsection{\@startsection{subsection}{3}
  \z@{.5\linespacing\@plus.7\linespacing}{-.5em}%
  {\normalfont\bfseries}}
  
  \makeatother
  
  \makeatletter
\newcommand{\crefnames}[3]{%
  \@for\next:=#1\do{%
    \expandafter\crefname\expandafter{\next}{#2}{#3}%
  }%
}
\makeatother

\newcommand{\nref}[1]{\nameref{#1} \ref{#1}}
\crefnames{part,chapter,section}{\S\!}{\S\S\!}

\newtheorem{proposition-definition}[theorem]{Proposition-Definition}

\theoremstyle{definition}

\makeatletter
\newtheorem*{rep@theorem}{\rep@title}
\newcommand{\newreptheorem}[2]{%
	\newenvironment{rep#1}[1]{%
		\def\rep@title{#2 \ref{##1}}%
		\begin{rep@theorem}}%
		{\end{rep@theorem}}}
\makeatother

\newreptheorem{theorem}{Theorem}
\newreptheorem{lemma}{Lemma}
\newreptheorem{corollary}{Corollary}
\newreptheorem{proposition}{Proposition}
\newreptheorem{setup}{Setup}

\setlist{leftmargin=15pt,labelindent=15pt}
\setlist[enumerate]{wide=\parindent, leftmargin=15pt, labelwidth=15pt, align=left,label=(\roman*)}

\counterwithin*{footnote}{subsection}

\mathchardef\mhyphen="2D 
\newlist{subenumerate}{enumerate}{1}
\setlist[subenumerate,1]{label=(\arabic*)}

\let\originalleft\left
\let\originalright\right
\renewcommand{\left}{\mathopen{}\mathclose\bgroup\originalleft}
\renewcommand{\right}{\aftergroup\egroup\originalright}

\mathchardef\mhyphen="2D 

\newcommand\defeq{\stackrel{\mathrm{\mbox{\scriptsize{def}}}}{\,=\,}}
\newcommand{\sth}{\,\,\vert\,\,}

\newcommand{\Z}{\mathbb{Z}}
\newcommand{\Q}{\mathbb{Q}}
\newcommand{\R}{\mathbb{R}}
\newcommand{\C}{\mathbb{C}}

\newcommand{\field}{\mathit{k}}
\newcommand{\modstack}{\mathscr{M}}

\newcommand{\calW}{\mathcal{W}}

\newcommand{\bfg}{\mathbf{g}}

\DeclareMathOperator{\fgMod}{\mathbf{mod}}
\DeclareMathOperator{\add}{\mathbf{add}}
\DeclareMathOperator{\filt}{\mathbf{filt}}
\DeclareMathOperator{\fac}{\mathbf{fac}}

\DeclareMathOperator{\Jasso}{\mathbf{J}}

\DeclareMathOperator{\taurigidpair}{\tau\mhyphen rigid{\,}pair}

\DeclareMathOperator{\ftors}{f\mhyphen tors}

\DeclareMathOperator{\Hom}{Hom}

\newcommand{\tcmc}[1]{\mathfrak{W}_{#1}}
\newcommand{\picspace}[1]{\mathfrak{X}_{#1}}
\newcommand{\picgroup}[1]{G_{#1}}

\newcommand{\VerdierProd}{\ast}
\newcommand{\HallProd}{\bullet}

\newcommand{\WL}{\mathbf{WL}}

\definecolor{friendly_deepblue}{RGB}{51,34,136}
\definecolor{friendly_blue}{RGB}{136,204,238}
\definecolor{friendly_beige}{RGB}{221,204,119}
\definecolor{friendly_red}{RGB}{136,34,85}
\definecolor{friendly_green}{RGB}{68,170,153}

\title[Picture groups and green sequences from Ringel--Hall algebras]{Picture groups and green sequences \\ from the perspective of Ringel--Hall algebras}

\author[Erlend D. Børve]{Erlend D. Børve}

\address{Institut for Matematik, Aarhus Universitet, Ny Munkegade 118, 8000 Aarhus C, Denmark}
\email{erlend.d.borve@math.au.dk}
\thanks{The author gratefully acknowledges support from
	the French ANR grant CHARMS (ANR-19-CE40-0017-02),
	the Deutsche Forschunggemeinschaft (DFG, German Research Foundation) -- Project ID 281071066 -- TRR 191, 
	the NAWI Graz Postdoc fellowship, 
	the Deutsche Forschungsgemeinschaft (DFG, German Research Foundation) -- Projektnummer 496500943, 
	and the Aarhus University Research Foundation -- Grant number AUFF-E-2024-9-43.}

\begin{document}
\begin{abstract}
Let $\field$ be an algebraically closed field and let $\Lambda$ be a finite-dimensional associative $\field$-algebra. 
We apply Joyce and Brideland's notion of Ringel--Hall algebra to prove results about picture spaces and picture groups. For instance, we construct a \textit{faithful group functor for $\Lambda$}, i.e. a faithful functor from the $\tau$-cluster morphism category of $\Lambda$ into a groupoid. Consequently, by results of Hanson--Igusa, the picture space of $\Lambda$ is locally CAT(0) provided that the $\tau$-cluster morphism category of $\Lambda$ admits compatibility of last factors. We also show that green sequences are in bijection with certain positive expressions in the picture group, which generalizes a result of Igusa--Todorov beyond hereditary $\field$-algebras.
\end{abstract}
	
	\maketitle
	\tableofcontents	
	
\section{Introduction}

\subsection{} \textit{Cluster algebras} \cite{FZ02i,FZ02ii,BFZ05,FZ07} have found a wide range of applications since the turn of the millennium. One instance is in homological mirror symmetry, where the duality of cluster varieties enables a reformulation of the central conjecture of the field \cite{FG06,HK18}. 
The notion of \textit{scattering diagram} later arose, and was used to settle several conjectures for cluster algebras, concerning positivity and sign coherence \cite{GHKK18}. One can generalize the construction of scattering diagrams to the setting of finite-dimensional algebras (and more generallly, to Artin algebras \cite{Tre26}), and recover them from the \textit{Ringel--Hall algebra} of the moduli stack of representations \cite{KS08,Bri17}. 

For a {representation-finite} hereditary algebra $\Lambda$, the underlying polyhedral fan of the scattering diagram (known as the \textit{wall-and-chamber structure} of $\Lambda$ \cite{BST19}) is equivalent to the \textit{semi-invariant picture} of $\Lambda$ \cite{IOTW15}. The semi-invariant picture was later used to define the \textit{picture group} of $\Lambda$ in terms of generators and relations \cite{ITW16}. The \textit{picture space} of $\Lambda$ is a finite CW-complex $\picspace{\Lambda}$ which is a $\mathrm{K}(\pi,1)$ for the picture group \cite{ITW16}. This is to say that the fundamental group $\pi_{1}(\picspace{\Lambda})$ is isomorphic to the picture group of $\Lambda$ and that the homotopy groups $\pi_{n}(\picspace{\Lambda})$ vanish for $n\geq 2$.

To systemitize the study of picture groups, Igusa--Todorov defined the cluster morphism category \cite{IT17}, a category whose classifying space is homeomorphic to the picture space and whose fundamental group is isomorphic to the picture group. Buan--Marsh \cite{BM18w} and Buan--Hanson \cite{BH21} later introduced the notion of \textit{$\tau$-cluster morphism category}, which enables the study of picture groups over any finite-dimensional algebra (further generalizations and equivalent definitions appeared later \cite{Bor21,STTW23,Kai23,Bor24,Kai24,BK26}). For an arbitrary finite-dimensional algebra, one can define the \textit{picture space} of $\Lambda$ as the classifying space of the $\tau$-cluster morphism category and the \textit{picture group} as a subgroup of the fundamental group of the picture space (see \Cref{picdefs}). It is not clear whether the picture space remains a $\mathrm{K}(\pi,1)$ for the picture group in general, nor does a counter-example seem to be known.

\subsection{} The purpose of this note is to introduce algebro-geometric techniques to the study of picture groups and picture spaces. The relationship between picture groups and Ringel--Hall algebras has long been known to experts, but a systematic treatment seems missing.\footnote{The connection between picture groups and scattering diagrams appears in recent work of Treffinger \cite{Tre26}.} Working over an algebraically closed field, we will map the picture group into the Ringel--Hall algebra of stack functions, so that we attain a perspective which is invisible from the presentation by generators and relations.
In the appendix of previous work with M. Kaipel, E. J. Hanson achieves the same for finite fields using classical Ringel--Hall algebras \cite[Appendix A]{BHK25}.

\subsection{Summary of results} 

In \Cref{sec:pic}, we recall the definitions of $\tau$-cluster morphism categories, in the sense of Buan--Hanson \cite{BH21}. For a finite-dimensional algebra $\Lambda$ over a field $\field$, let $\tcmc{\Lambda}$ denote the $\tau$-cluster morphism category of $\Lambda$. In  \Cref{picdefs}, the \textit{picture space} $\picspace{\Lambda}$ of $\Lambda$ is defined as a the classifying space of $\tcmc{\Lambda}$, and the \textit{picture group} $\picgroup{\Lambda}$ of $\Lambda$ will be defined as the subgroup of the fundamental group of $\picspace{\Lambda}$. The necessary background on Ringel--Hall algebras is recalled in \Cref{sec:RingelHall}, where we review work of Joyce \cite{Joy07} and Bridgeland \cite{Bri17}. Of particular importance is the pro-unipotent group $\widehat{G}$, which can be identified with $1_{\mathbf{0}}+\widehat{H}_{>0}(\Lambda)$, where $\widehat{H}_{>0}(\Lambda)$ is the positive graded part of the completed Ringel--Hall algebra $\widehat{H}(\Lambda)$ and $1_{\mathbf{0}}\in\widehat{H}(\Lambda)$ is the identity element.

Our first main result is the following:\footnote{\nref{thm:2main} was proven independently by Treffinger \cite[§6]{Tre26}.}

\begin{reptheorem}{thm:2main}
Let $\field$ be an algebraically closed field and let $\Lambda$ be a finite-dimensional $\field$-algebra. 
\begin{enumerate}
\item\label{thm:2main1} Regard the group $\widehat{G}$ as a groupoid with a single object. We may define a faithful functor $\eta\colon \tcmc{\Lambda} \to \widehat{G}$. 
\item\label{thm:2main2} The generators $g_B$ of the picture group $\picgroup{\Lambda}$, where $B$ runs over the f-bricks in $\fgMod(\Lambda)$, are pair-wise distinct.
\end{enumerate}
\end{reptheorem}

One of the significant consequences of \nref{thm:2main}\ref{thm:2main1} concerns the homotopy type of picture spaces. The $\tau$-cluster morphism category $\tcmc{\Lambda}$ is a cubical category, in the sense of Igusa \cite[§3.1]{Igu14}. Using work of Gromov \cite{Gro87}, Igusa provides sufficient conditions for a cubical category to have a classfying space which is locally $\mathrm{CAT}(0)$ \cite[Proposition 3.5]{Igu14}. We recite these conditions in \Cref{Igu14.3.4}. The first condition \ref{I1} always holds for $\tau$-cluster morphism categories, and is easy to verify using $\tau$-tilting theory. The third condition, namely \ref{I3}, has been shown in \nref{thm:2main}\ref{thm:2main1} to hold over any algebraically closed field (with M. Kaipel, we generalize to perfect fields \cite[Corollary 4.13]{BHK25}). We are left with \ref{I2}, which is known as \textit{compatibility of last factors}.

\begin{repcorollary}{cor:complastK_alg}
Let $\field$ be an algebraically closed field and let $\Lambda$ be a finite-dimensional $\field$-algebra. 
If compatibility of last factors holds for $\Lambda$, then the picture space $\picspace{\Lambda}$ is locally $\mathrm{CAT}(0)$. In particular, the condition \ref{I2} is sufficient for $\picspace{\Lambda}$ to be an Eilenberg--MacLane space of type $\mathrm{K}(\pi,1)$.
 \end{repcorollary}

In the representation-finite hereditary case, Igusa--Todorov provide a bijeciton between positive expressions in the picture group and green sequences \cite{IT21}. In \Cref{sec:green}, we generalize Igusa--Todorov's result, showing that it holds for any finite-dimensional algebra over an algebraically closed field.

\begin{reptheorem}{thm:positive_expr}
Let $\field$ be an algebraically closed field and let $\Lambda$ be a finite-dimensional $\field$-algebra. Given a green sequence
 \begin{equation*}
	\mathbf{0} =  \mathcal{T}_0 \lessdot \mathcal{T}_{1} \lessdot \cdots \lessdot \mathcal{T}_{\ell-1} \lessdot \mathcal{T}_{\ell} = \mathcal{T}
\end{equation*}
in $\fgMod(\Lambda)$, with brick labels reading $(B_{1}, \dots , B_{\ell-1},B_{\ell})$ from left to right, let $g_i$ denote the generator of the picture group $\picgroup{\Lambda}$ provided by the brick $B_i$ (see \eqref{eq:picgroup_pres_genloop}). 
Mapping the green sequence in \eqref{eq:greenseqThm} to the word $g_{1}\cdots g_{\ell-1} g_{\ell}$ of generators of $\picgroup{\Lambda}$ defines a bijection from the first to the second of the following sets:
\begin{enumerate}
	\item green sequences in $\fgMod(\Lambda)$ of the form displayed above,
	\item positive expressions in the picture group $\picgroup{\Lambda}$ (i.e. a word of generators) whose product is equal to $g_{\mathcal{T}}$ as group elements, where $g_{\mathcal{T}}$ is the product $g_{1}\cdots g_{\ell-1} g_{\ell}$ in $\picgroup{\Lambda}$.
\end{enumerate}
In particular, we have a bijection from the first to the second of the following sets:
\begin{enumerate}
	\item maximal green sequences in $\fgMod(\Lambda)$,
	\item positive expressions in the picture group $\picgroup{\Lambda}$ whose product is equal to $g_{\fgMod(\Lambda)}$.
\end{enumerate}
\end{reptheorem}

\subsection{Notation and conventions} 
Throughout the entire text, we fix a field $\field$ and a finite-dimensional associative $\field$-algebra $\Lambda$. By \textit{$\Lambda$-module}, we will always mean finitely generated right $\Lambda$-module. Let $\fgMod(\Lambda)$ denote the category of finitely generated right $\Lambda$-modules. 
An \textit{additive subcategory} of an additive category $\mathcal{A}$ is full subcategory of $\mathcal{A}$ which is closed under isomorphisms, direct sums and direct summands. If $\mathcal{X}$ is a full subcategory of $\mathcal{A}$, let $\add(\mathcal{X})$ denote the smallest additive subcategory of $\mathcal{A}$ containing $\mathcal{X}$. If $\mathcal{X}$ consists of a single object $X$, we write $\add(X)$ when we mean $\add(\mathcal{X})$.

\subsection{Acknowledgements} 
The author thanks Tom Bridgeland, Carlo Klapproth, Monica Garcia, Maximilian Kaipel, Bernhard Keller, and Lang Mou for insightful discussions. Eric J. Hanson's comments on a previous draft are also gratefully acknowledged.

\section{$\tau$-cluster morphism categories and picture groups}\label{sec:pic}
\setcounter{subsection}{-1}

\subsection{} In the literature, the notion of $\tau$-cluster morphism categories has been defined in various levels of generality \cite{Igu14,CN17,IT17,BM18w,BH21,Bor21,STTW23,Kai23,Bor24,Kai24,Non26,BK26}. This section reviews $\tau$-cluster morphism categories and related notions. Here, the field $\field$ is arbitrary. In particular, it need not be algebraically closed.

\subsection{}\label{1.1}
Let $\mathcal{X}$ and $\mathcal{Y}$ be additive subcategories of $\fgMod(\Lambda)$. The \textit{Verdier product} of $\mathcal{X}$ and $\mathcal{Y}$ is the full subcategory of $\fgMod(\Lambda)$ defined by 
\begin{equation*}
 \mathcal{X} \VerdierProd \mathcal{Y} 
 \defeq \{ M\in \fgMod(\Lambda) \sth \text{ $\exists$ a short exact sequence } 
 \begin{tikzcd}[column sep=1.5em] M_{\mathcal{X}}\arrow[r,tail] & M \arrow[r,two heads] & M_{\mathcal{Y}}, \end{tikzcd}
 \text{ with $M_{\mathcal{X}} \in {\mathcal{X}}$ and $M_{\mathcal{Y}} \in {\mathcal{Y}}$}   \}.
\end{equation*}
Using the Third Isomorphism Theorem, it can be shown that $\VerdierProd$ is associative. We say that an additive subcategory $\mathcal{X}\subseteq \fgMod(\Lambda)$ is \textit{extension-closed} (or \textit{closed under extensions}) of $\mathcal{X}\VerdierProd \mathcal{X} \subseteq \mathcal{X}$.

Let $\mathcal{X}$ be an additive subcategory of $\fgMod(\Lambda)$. 
{A \textit{right $\mathcal{X}$-approximation} (resp. \textit{left $\mathcal{X}$-approximation}) of a $\Lambda$-module $M$ is a $\Lambda$-homomorphism $M_{\mathcal{X}}\to^{p}M$ (resp. $M\to^{i} M_{\mathcal{X}}$) in $\fgMod(\Lambda)$ such that $M_{\mathcal{X}}\in \mathcal{X}$ and the induced $k$-homomorphism $\Hom_{\Lambda}(X,M_{\mathcal{X}}) \to^{p\circ -}\Hom_{\Lambda}(X,M)$ (resp. $\Hom_{\Lambda}(M_{\mathcal{X}},X)\to^{- \circ i}\Hom_{\Lambda}(M,X) $) is surjective whenever $X\in \mathcal{X}$. We say that $\mathcal{X}$ is a \textit{contravariantly finite} (resp. \textit{covariantly finite}, resp. \textit{functorially finite}) if every object in $\fgMod(\Lambda)$ admits a right (resp. left, resp. right and left) $\mathcal{X}$-approximation. 

A \textit{torsion class} (resp. \textit{torsion-free class}) of $\fgMod(\Lambda)$ is an additive subcategory $\mathcal{T}\subseteq \fgMod(\Lambda)$ which is closed under extensions and factor modules (resp. extensions and submodules) \cite{Dic66}. Given a torsion class $\mathcal{T}$, its \textit{right perpendicular category}, namely
\begin{equation*}
 \mathcal{T}^{\perp} \defeq \{M\in \fgMod(\Lambda) \sth \Hom_{\Lambda}(X,M)=0 \quad \forall X\in\mathcal{T} \}
\end{equation*}
is a torsion-free class. We have that $\mathcal{T}$ is functorially finite precisely when $\mathcal{T}^{\perp}$ is functorially finite \cite{Sma84}.
A pair of the form $(\mathcal{T},\mathcal{T}^{\perp})$, where $\mathcal{T}$ is a torsion class, is called a \textit{torsion pair}. 
For a torsion class $\mathcal{T}$ and a $\Lambda$-module $M$, there is a unique short exact sequence in $\fgMod(\Lambda)$ up to isomorphism
\begin{equation}\label{eq:canseq}
	\begin{tikzcd}
		M_{\mathcal{T}}\arrow[r,tail] & M \arrow[r,two heads] & M_{\mathcal{T}^{\perp}},
	\end{tikzcd}
\end{equation}
with the property that $M_{\mathcal{T}} \in \mathcal{T}$ and $M_{\mathcal{T}^{\perp}}\in \mathcal{T}^{\perp}$. In particular, we have that $\mathcal{T}\VerdierProd\mathcal{T}^{\perp} = \fgMod(\Lambda)$. The assignment $M\mapsto M_{\mathcal{T}}$ (resp. $M\mapsto M_{\mathcal{T}^{\perp}}$) determines an additive functor $\fgMod(\Lambda) \to \mathcal{T}$ (resp. $\fgMod(\Lambda) \to \mathcal{T}^{\perp}$).

\subsection{}\label{1.2}
A \textit{wide subcategory} of $\fgMod(\Lambda)$ is an extension-closed abelian subcategory of $\fgMod(\Lambda)$. Equivalently, \textit{wide subcategory} of $\fgMod(\Lambda)$ is an additive subcategory of $\fgMod(\Lambda)$ which is closed under extensions, kernels, and cokernels. Let $\mathrm{wide}(\Lambda)$ denote the poset of wide subcategories of $\fgMod(\Lambda)$.

For an additive subcategory $\mathcal{X}$ of $\fgMod(\Lambda)$, let $\filt(\mathcal{X})$ denote the smallest extension-closed additive subcategory of $\fgMod(\Lambda)$ containing $\mathcal{X}$, let $\fac(\mathcal{X})$ denote the closure of $\mathcal{X}$ under factor modules 
We write $\filt(X)$ and $\fac(X)$ whenever $\mathcal{X}=\add(X)$ for some $\Lambda$-module $X$.

One can show that a wide subcategory $\calW$ of $\fgMod(\Lambda)$ is functorially finite precisely when there exists a finite-dimensional $\field$-algebra $\Gamma$ and an exact equivalence $\calW \to^{\simeq} \fgMod(\Gamma)$ \cite[Proposition 4.12]{Eno22}. Since a functorially finite wide subcategory $\mathcal{W}$ of $\fgMod(\Lambda)$ is equivalent to the module category of a finite-dimensional $\field$-algebra $\Gamma$, one can implicity fix $\Gamma$ and define all operations in above for $\calW$. We will adapt the notation by adding $\calW$ as a subscript, whence we write $\filt_{\calW}$ and $\fac_{\calW}$.

\subsection{}\label{ttilstuff}
Let $\tau$ denote the Auslander--Reiten translation in $\fgMod(\Lambda)$. Following Adachi--Iyama--Reiten \cite{AIR14}, we say that a $\Lambda$-module $M$ is \textit{$\tau$-rigid} if $\Hom_{\Lambda}(M,\tau M)=0$. A pair $(M,P)$ of $\Lambda$-modules is a \textit{$\tau$-rigid pair} if $M$ is $\tau$-rigid and $P$ is a projective $\Lambda$-module satisfying $\Hom_{\Lambda}(P,M)=0$. A $\tau$-rigid pair $(M,P)$ is a \textit{direct summand} of the $\tau$-rigid pair $(N,Q)$ if $M$ is a direct summand of $N$ and $P$ is a direct summand of $Q$. 
A $\tau$-rigid pair $(M,P)$ is \textit{basic} if both $M$ and $P$ are basic $\Lambda$-modules, and \textit{indecomposable} if $(0,0)$ is its only proper direct summand. A \textit{support $\tau$-tilting pair} is a $\tau$-rigid pair with a maximal number of direct summands up to isomorphism.

We denote the set of basic $\tau$-rigid pairs in $\fgMod(\Lambda)$ by $\taurigidpair(\Lambda)$, and those having a $\tau$-rigid pair $(M,P)$ as a direct summand by $\taurigidpair_{(M,P)}(\Lambda)$. 
Since a functorially finite wide subcategory $\calW$ of $\fgMod(\Lambda)$ is equivalent to the module category of a finite-dimensional $\field$-algebra, we may define the notion of $\tau$-rigid pairs in $\calW$. The sets $\taurigidpair(\calW)$ and $\taurigidpair_{(M,P)}(\calW)$ are then defined in the same spirit as above.

If $M$ is a $\tau$-rigid $\Lambda$-module, then $\fac(M)$ is a functorially finite torsion class in $\fgMod(\Lambda)$ \cite[Theorem 5.10]{AS81}, and $M^{\perp}$ is the corresponding torsion-free class. In the following, the canonical short exact sequences for the torsion pair $(\fac(M),M^{\perp})$ (see \eqref{eq:canseq}) will be denoted
\begin{equation*}
\begin{tikzcd}
 i_M(X)\arrow[r,tail]& X \arrow[r,two heads] & p_M(X),
\end{tikzcd}
\end{equation*}
which is to say that $i_M\colon \fgMod(\Lambda)\to \fac(M)$ and $p_M\colon \fgMod(\Lambda)\to M^{\perp}$ denote the functorial assignments into the torsion class and the torsion-free class, respectively. More generally, if $M$ is an object in a functorially finite wide subcategory $\mathcal{W}$ of $\fgMod(\Lambda)$, we also denote the functorial assignments by
$i_M\colon \mathcal{W} \to \fac_{\mathcal{W}}(M)$ and $p_M\colon \mathcal{W} \to M^{\perp}\cap \mathcal{W}$, respectively.

\subsection{}\label{ttilstuff'}
Let $\mathcal{W}$ be a functorially finite wide subcategory of $\fgMod(\Lambda)$ and let $(M,P)$ be a $\tau$-rigid pair in $\mathcal{W}$. The \textit{$\tau$-perpendicular category of} $(M,P)$ in $\calW$ is defined by $\Jasso_{\mathcal{W}}(M,P)\defeq M^{\perp} \cap {^{\perp}(\tau_{\mathcal{W}}M})\cap P^{\perp}$, where $\tau_{\mathcal{W}}$ denotes the Auslander--Reiten translation in $\calW$. This is a wide subcategory of $\calW$ \cite[Proposition 3.6]{Jas15} \cite[Theorem 4.12(a)]{DIRRT17} \cite[Corollary 3.25]{BST19}. The wide subcategories of $\calW$ arising in this way are the \textit{$\tau$-perpendicular subcategories} of $\calW$. When $\calW=\fgMod(\Lambda)$, we remove the subscript and write $\Jasso(M,P)$ for the {$\tau$-perpendicular category of} $(M,P)$.
For all $\tau$-rigid pairs $(M,P)$ in $\mathcal{W}$, there is a bijection \cite[§3]{BM18w}
\begin{equation}\label{BM-bij}
	E_{(M,P)}\colon \taurigidpair_{(M,P)}(\mathcal{W}) \to \taurigidpair(\Jasso_{\mathcal{W}}(M,P)).
\end{equation}
This bijection is defined in a case-by-case manner, but we do not give the full definition here.

\subsection{Lemma}\label{lem:facMM'}
	Let $\mathcal{W}$ be a functorially finite wide subcategory of $\fgMod(\Lambda)$, let $(M,P)$ be a $\tau$-rigid pair in $\mathcal{W}$, and let $\mathcal{V}$ denote $\Jasso_{\mathcal{W}}(M,P)$.
	If $E_{(M,P)}^{-1}(N,Q) = (M\oplus M',P\oplus P')$, then $\fac_{\mathcal{W}}(M\oplus M') = \fac_{\mathcal{W}}(M) \VerdierProd \fac_{\mathcal{V}}(N)$.
\begin{proof}
	In all cases defining the map $E_{(M,P)}$, we have the following: if $E_{(M,P)}^{-1}(N,Q) = (M\oplus M',P\oplus P')$, then $N=p_{M}(M\oplus M')$. Since $\fac_{\mathcal{W}}(M\oplus M')$ is a torsion class in $\mathcal{W}$, we deduce that \[ \fac_{\mathcal{W}}(M\oplus M') \supseteq \fac_{\mathcal{W}}(M) \VerdierProd \fac_{\mathcal{V}}(N).\] For the reverse inclusion, consider the solid part of the following diagram
	\begin{equation*}
		\begin{tikzcd}
			i_{M}(M\oplus M')^{\oplus m} \arrow[r,tail] & (M\oplus M')^{\oplus m}  \arrow[r,two heads] \arrow[d,two heads]& N^{\oplus m}\arrow[d,two heads,dashed] \\
			i_{M}(X) \arrow[r,tail] & X  \arrow[r,two heads] & p_{M}(X)
		\end{tikzcd}
	\end{equation*}
	where $X\in \fac_{\mathcal{W}}(M\oplus M')$ and the rows are canonical short exact sequences for the torsion pair $(\fac_{\mathcal{W}}(M),M^{\perp}\cap \mathcal{W})$ in $\mathcal{W}$. Having produced the bottom row above, it now suffices to show that $p_{M}(X)\in \mathcal{V}$. Indeed, it is in $M^{\perp}\cap \mathcal{W}$ by construction, and by the exactness properties of Hom-functors, the existence of the dashed epimorphism above is used to show that $p_{M}(X)\in {^{\perp}(\tau_{\mathcal{W}}M)}\cap P^{\perp}\cap \mathcal{W}$.
\end{proof}

\subsection{} 
A $\Lambda$-module $B$ is a \textit{brick} if the endomorphism $\field$-algebra $\mathrm{End}_{\Lambda}(B)$ is a division $\field$-algebra. If $B$ is contained in a functorially finite wide subcategory $\mathcal{W}$ of $\fgMod(\Lambda)$, we may say that $B$ is  a \textit{brick in $\mathcal{W}$}.
We say that $\Lambda$ is \textit{brick-finite} if admits only finitely many bricks up to isomorphism.
Given a brick $B$ in $\mathcal{W}\subseteq\fgMod(\Lambda)$, we call $B$ an \textit{f-brick} if $\filt_{\mathcal{W}}(B)$ is a $\tau$-perpendicular subcategory of $\mathcal{W}$. 

\subsection{}\label{picdefs}
Buan--Hanson \cite{BH21} define the \textit{$\tau$-cluster morphism category} $\tcmc{\Lambda}$ of $\Lambda$ as follows. The objects are the $\tau$-perpendicular subcategories of $\fgMod(\Lambda)$ and a morphism $\mathcal{W} \to^{(M,P)} \mathcal{V}$ is given by a basic $\tau$-rigid pair $(M,P)$ in $\mathcal{W}$ such that $\Jasso_{\mathcal{W}}(M,P)=\mathcal{V}$. 
Given two composable morphisms
\begin{equation*}
\begin{tikzcd}
	\mathcal{W} \arrow[r,"{(M,P)}"] & \mathcal{V}\arrow[r,"{(N,Q)}"] & \mathcal{V}',
\end{tikzcd}
\end{equation*}
their composite is given by 
\begin{tikzcd}[column sep=5em]
	\mathcal{W} \arrow[r,"{E_{(M,P)}^{-1}(N,Q)}"] & \mathcal{V}',
\end{tikzcd}
where $E$ is as in \eqref{BM-bij}.

The \textit{picture space} of $\Lambda$ is the classifying space of the $\tau$-cluster morphism category $\tcmc{\Lambda}$. We denote it by $\picspace{\Lambda}$. For any object $\calW$ in $\tcmc{\Lambda}$, we have a morphism $\calW \to \mathbf{0}$, where $\mathbf{0}$ denotes the wide subcategory of $\fgMod(\Lambda)$ consisting of zero objects.
This makes $\mathbf{0}$ a weakly terminal object in $\tcmc{\Lambda}$, whence the picture space $\picspace{\Lambda}$ becomes a path-connected topological space. 
The \textit{picture group} of $\Lambda$ is denoted $\picgroup{\Lambda}$, and defined as the subgroup of the fundamental group $\pi_1 (\tcmc{\Lambda},\mathbf{0})$ generated by loops of the form 
\begin{equation}\label{eq:picgroup_pres_genloop}
\begin{tikzcd}[column sep=3em]
			 \mathbf{0} & \arrow[l,"{(0,T)}"'] \filt(B) \arrow[r,"{(T,0)}"] & \mathbf{0},
		\end{tikzcd}
\end{equation}
where $B$ is an f-brick in $\fgMod(\Lambda)$ and $T$ is a projective cover of $B$ in $\filt(B)$ (since $\filt(B)$ is equivalent to the module category of a local $k$-algebra, we have that $T$ is the unique $\tau$-rigid object in $\filt(B)$). The generator of the picture group $\picgroup{\Lambda}$ displayed in \eqref{eq:picgroup_pres_genloop} will be denoted $g_B$. Our definition is more general than that Hanson--Igusa, but we recover their definition (up to isomorphism) when $\Lambda$ is brick-finite \cite[Theorem 4.10]{HI21}.

\subsection{Examples} 

\subsubsection{} Suppose that $\Lambda$ is a local finite-dimensional $\field$-algebra. There are only two non-zero basic $\tau$-rigid pairs in $\fgMod(\Lambda)$, namely $(\Lambda,0)$ and $(0,\Lambda)$. The $\tau$-cluster morphism category of $\Lambda$ becomes
\begin{equation*}
	\begin{tikzcd}[column sep=5em,
    execute at end picture={
    \scoped[on background layer]
    \fill[color=friendly_green, opacity=0.5] (SP.north west) -- (SP.north east)  -- (SP.south east) -- (SP.south west) -- cycle;},
     execute at end picture={
    \scoped[on background layer]
    \fill[color=friendly_green, opacity=0.5] (P.north west) -- (P.north east)  -- (P.south east) -- (P.south west) -- cycle;}]
	|[alias=SP]| \mathbf{0} &  |[alias=o]|\fgMod(\Lambda) \arrow[r, "{\tiny (\Lambda,0)}"] \arrow[l, "{\tiny (0,\Lambda)}"']  & |[alias=P]| \mathbf{0}
	\end{tikzcd}
\end{equation*}
where the two occurrences of $\mathbf{0}$ (shaded in \textcolor{friendly_green}{green}) are to be identified. The picture space of $\Lambda$ is homeomorphic to a circle, and the picture group of $\Lambda$ is free abelian with one generator.

\subsubsection{} Suppose that $\Lambda$ is the path $\field$-algebra of the quiver  $Q = \begin{tikzcd} 1\arrow[r,"\alpha"] & 2. \end{tikzcd}$
The Auslander--Reiten quiver of the finite-dimensional $\field$-algebra $\Lambda $ is then
\begin{equation}\label{eg:ARquiv}
\begin{tikzcd}
               & P_2 \arrow[rd] &                \\
P_1 \arrow[ru] &                & S_2 \arrow[ll,dashed,"\tau"']
\end{tikzcd}
\end{equation}
where $P_1$ is simple and projective, $P_2$ is non-simple and projective-injective, and $S_2$ is simple and injective. The $\field$-algebra $\Lambda$ is of finite representation type, so it is {brick-finite} in particular. All three indecomposable $\Lambda$-modules are $\tau$-rigid. 
The $\tau$-cluster morphism category of $\Lambda$ is the following
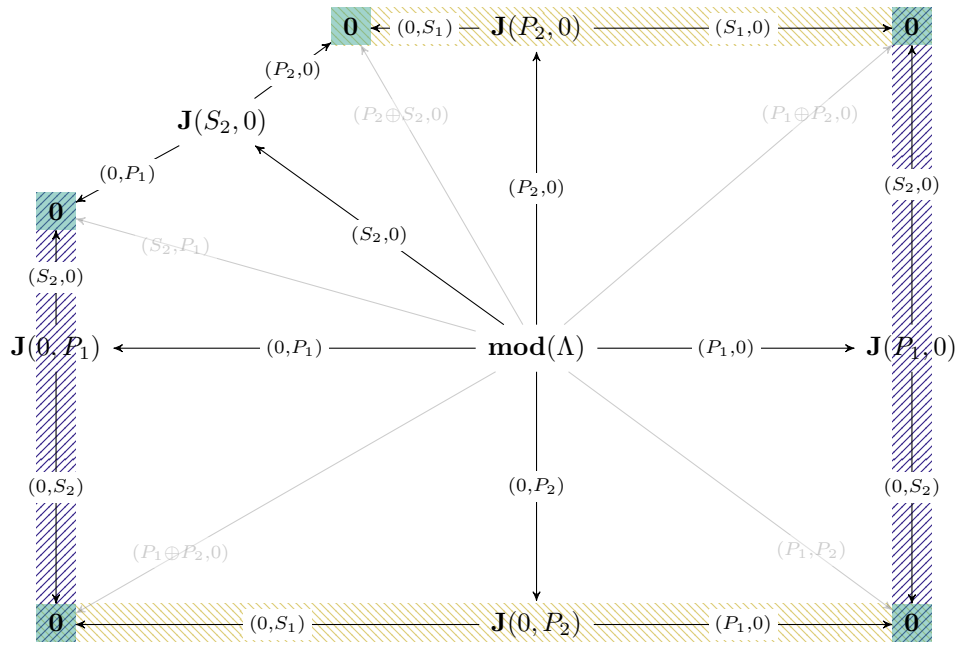
\begin{equation*}
    \begin{tikzcd}[column sep=2em,
    execute at end picture={
    \scoped[on background layer]
    execute at end picture={
    \scoped[on background layer]
    \fill[pattern=north west lines, pattern color=friendly_beige] (p2s2.north west) -- (p1p2.north east)  -- (p1p2.south east) -- (p2s2.south west) -- cycle;},
    execute at end picture={
    \scoped[on background layer]
    \fill[pattern=north west lines, pattern color=friendly_beige] (Sp1Sp2.north west) -- (s1Sp2.north east)  -- (s1Sp2.south east) -- (Sp1Sp2.south west) -- cycle;},
    execute at end picture={
    \scoped[on background layer]
    \fill[pattern=north east lines, pattern color=friendly_deepblue] (s2Sp1.north west) -- (Sp1Sp2.south west)  -- (Sp1Sp2.south east) -- (s2Sp1.north east) -- cycle;},
    execute at end picture={
    \scoped[on background layer]
    \fill[pattern=north east lines, pattern color=friendly_deepblue] (p1p2.north west) -- (s1Sp2.south west)  -- (s1Sp2.south east) -- (p1p2.north east) -- cycle;},
    \fill[color=friendly_green, opacity=0.5] (s1Sp2.north west) -- (s1Sp2.north east)  -- (s1Sp2.south east) -- (s1Sp2.south west) -- cycle;},
    execute at end picture={
    \scoped[on background layer]
    \fill[color=friendly_green, opacity=0.5] (Sp1Sp2.north west) -- (Sp1Sp2.north east)  -- (Sp1Sp2.south east) -- (Sp1Sp2.south west) -- cycle;},
    execute at end picture={
    \scoped[on background layer]
    \fill[color=friendly_green, opacity=0.5] (s2Sp1.north west) -- (s2Sp1.north east)  -- (s2Sp1.south east) -- (s2Sp1.south west) -- cycle;},
    execute at end picture={
    \scoped[on background layer]
    \fill[color=friendly_green, opacity=0.5] (p2s2.north west) -- (p2s2.north east)  -- (p2s2.south east) -- (p2s2.south west) -- cycle;},
    execute at end picture={
    \scoped[on background layer]
    \fill[color=friendly_green, opacity=0.5] (p1p2.north west) -- (p1p2.north east)  -- (p1p2.south east) -- (p1p2.south west) -- cycle;},
    ,
    ]
                                                                                                                                                                        &                                                                                                                                    & |[alias=p2s2]| \mathbf{0} &  & |[alias=p2]|\Jasso(P_2,0) \arrow[rrrr, "{\tiny (S_1,0)}"' description] \arrow[ll, "{\tiny (0,S_1)}" description]                                                                                                                                                                                                                                                                                                                                                                                                                                                                                                                                                         &  &                                &                                                                                                                                    & |[alias=p1p2]| \mathbf{0}                                                                                                      \\
                                                                                                                                                                        & |[alias=s2]|\Jasso(S_2,0) \arrow[ru, "{\tiny (P_2,0)}"' description] \arrow[ld, "{\tiny (0,P_1)}" description] &                         &  &                                                                                                                                                                                                                                                                                                                                                                                                                                                                                                                                                                                                                                                                                       &  &                                &                                                                                                                                    &                                                                                                                               \\
|[alias=s2Sp1]| \mathbf{0}                                                                                                                                          &                                                                                                                                    &                         &  &                                                                                                                                                                                                                                                                                                                                                                                                                                                                                                                                                                                                                                                                                       &  &                                &                                                                                                                                    &                                                                                                                               \\
                                                                                                                                                                        &                                                                                                                                    &                         &  &                                                                                                                                                                                                                                                                                                                                                                                                                                                                                                                                                                                                                                                                                       &  &                                &                                                                                                                                    &                                                                                                                               \\
|[alias=Sp1]|\Jasso(0,P_1) \arrow[uu, "{\tiny (S_2,0)}" description] \arrow[dddd, "{\tiny (0,S_2)}" description] &                                                                                                                                    &                         &  & |[alias=o]| \fgMod(\Lambda) \arrow[rrrr, "{\tiny{(P_1,0)}}" description] \arrow[dddd, "{\tiny {(0, P_2)}}" description] \arrow[llll, "{\tiny{(0,P_1)}}" description] \arrow[uuuu, "{\tiny {(P_2,0)}}" description] \arrow[llluuu, "{\tiny {(S_2,0)}}" description] \arrow[lluuuu, "{\tiny {(P_2\oplus S_2,0)}}" description,color=gray,opacity=0.4,near end] \arrow[rrrruuuu, "{\tiny{(P_1\oplus P_2,0)}}" description,color=gray,opacity=0.4,near end] \arrow[rrrrdddd, "{\tiny{(P_1,P_2)}}" description,color=gray,opacity=0.4,near end] \arrow[lllldddd, "{\tiny {(P_1\oplus P_2,0)}}" description,color=gray,opacity=0.4,near end] \arrow[lllluu, "{\tiny {(S_2, P_1)}}" description,color=gray,opacity=0.4,near end] &  &                                &                                                                                                                                    & |[alias=p1]| \Jasso(P_1,0) \arrow[uuuu, "{\tiny (S_2,0)}" description] \arrow[dddd, "{\tiny (0,S_2)}" description] \\
                                                                                                                                                                        &                                                                                                                                    &                         &  &                                                                                                                                                                                                                                                                                                                                                                                                                                                                                                                                                                                                                                                                                       &  &                                &                                                                                                                                    &                                                                                                                               \\
                                                                                                                                                                        &                                                                                                                                    &                         &  &                                                                                                                                                                                                                                                                                                                                                                                                                                                                                                                                                                                                                                                                                       &  &                                &                                                                                                                                    & |[alias=p1s1]| {}                                                                                               \\
                                                                                                                                                                        &                                                                                                                                    &                         &  &                                                                                                                                                                                                                                                                                                                                                                                                                                                                                                                                                                                                                                                                                       &  &                                & |[alias=s1]|{} &                                                                                                                               \\
|[alias=Sp1Sp2]|\mathbf{0}                                                                                                                                 &                                                                                                                                    &                         &  & |[alias=Sp2]|\Jasso(0,P_2) \arrow[llll, "{\tiny (0,S_1)}"' description] \arrow[rrrr, "{\tiny (P_1,0)}" description]                                                                                                                                                                                                                                                                                                                                                                                                                                                                                                                &  & {} &                                                                                                                                   &         |[alias=s1Sp2]| \mathbf{0}                                                                                                                       
\end{tikzcd}
\end{equation*}
The parts that have the same colour/shading are to be identified. The faint arrows are morphisms given by $\tau$-rigid pairs with two indecomposable direct summands, whereas the solid arrows are given by indecomposable $\tau$-rigid pairs.
The picture space of $\Lambda$ is homeomorphic to a punctured torus. Indeed, the puncture passes through the object $J(S_2,0)$.

\subsubsection{} Consider the quiver  $Q = \begin{tikzcd} 1\arrow[r,"\alpha", bend left=10,yshift=0.25em] & \arrow[l,"\beta", bend left=10,yshift=-0.25em] 2 \end{tikzcd}$ and the ideal $I\leq \field Q$ generated by the paths $\alpha\beta$ and $\beta\alpha$. 
The Auslander--Reiten quiver of the finite-dimensional $\field$-algebra $\Lambda = \field Q /I$ takes the form 
\begin{equation}\label{eg:ARquiv}
\begin{tikzcd}
& P_1 \arrow[rd] &                & P_2 \arrow[rd] &                & P_1 \arrow[rd] \\
\cdots\arrow[ru] &               & S_1 \arrow[ru] \arrow[ll,dashed,"\tau"'] &                & S_2 \arrow[ru]\arrow[ll,dashed,"\tau"'] &   & \cdots \arrow[ll,dashed,"\tau"'] 
\end{tikzcd}
\end{equation}
where the two vertices labelled $P_1$ are to be identified. All four indecomposable $\Lambda$-modules are $\tau$-rigid. 
The $\tau$-cluster morphism category of $\Lambda$ then becomes:
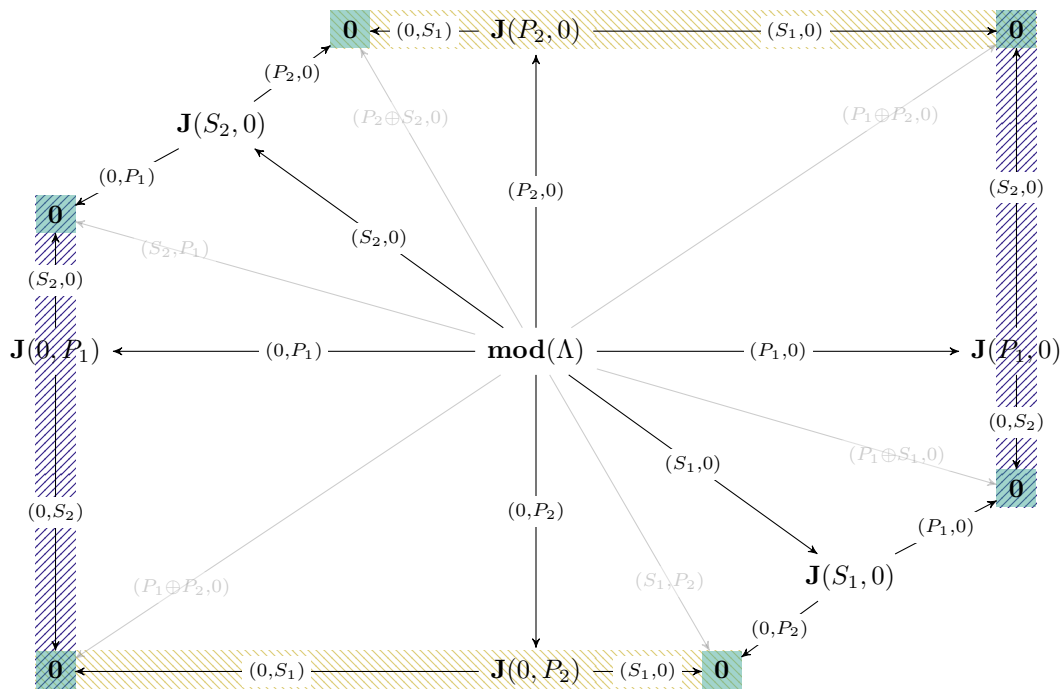
\begin{equation*}
    \begin{tikzcd}[column sep=2em,
    execute at end picture={
    \scoped[on background layer]
    execute at end picture={
    \scoped[on background layer]
    \fill[pattern=north west lines, pattern color=friendly_beige] (p2s2.north west) -- (p1p2.north east)  -- (p1p2.south east) -- (p2s2.south west) -- cycle;},
    execute at end picture={
    \scoped[on background layer]
    \fill[pattern=north west lines, pattern color=friendly_beige] (Sp1Sp2.north west) -- (s1Sp2.north east)  -- (s1Sp2.south east) -- (Sp1Sp2.south west) -- cycle;},
    execute at end picture={
    \scoped[on background layer]
    \fill[pattern=north east lines, pattern color=friendly_deepblue] (s2Sp1.north west) -- (Sp1Sp2.south west)  -- (Sp1Sp2.south east) -- (s2Sp1.north east) -- cycle;},
    execute at end picture={
    \scoped[on background layer]
    \fill[pattern=north east lines, pattern color=friendly_deepblue] (p1p2.north west) -- (p1s1.south west)  -- (p1s1.south east) -- (p1p2.north east) -- cycle;},
    \fill[color=friendly_green, opacity=0.5] (s1Sp2.north west) -- (s1Sp2.north east)  -- (s1Sp2.south east) -- (s1Sp2.south west) -- cycle;},
    execute at end picture={
    \scoped[on background layer]
    \fill[color=friendly_green, opacity=0.5] (Sp1Sp2.north west) -- (Sp1Sp2.north east)  -- (Sp1Sp2.south east) -- (Sp1Sp2.south west) -- cycle;},
    execute at end picture={
    \scoped[on background layer]
    \fill[color=friendly_green, opacity=0.5] (s2Sp1.north west) -- (s2Sp1.north east)  -- (s2Sp1.south east) -- (s2Sp1.south west) -- cycle;},
    execute at end picture={
    \scoped[on background layer]
    \fill[color=friendly_green, opacity=0.5] (p2s2.north west) -- (p2s2.north east)  -- (p2s2.south east) -- (p2s2.south west) -- cycle;},
    execute at end picture={
    \scoped[on background layer]
    \fill[color=friendly_green, opacity=0.5] (p1p2.north west) -- (p1p2.north east)  -- (p1p2.south east) -- (p1p2.south west) -- cycle;},
    execute at end picture={
    \scoped[on background layer]
    \fill[color=friendly_green, opacity=0.5] (p1s1.north west) -- (p1s1.north east)  -- (p1s1.south east) -- (p1s1.south west) -- cycle;},
    ]
                                                                                                                                                                        &                                                                                                                                    & |[alias=p2s2]| \mathbf{0} &  & |[alias=p2]|\Jasso(P_2,0) \arrow[rrrr, "{\tiny (S_1,0)}"' description] \arrow[ll, "{\tiny (0,S_1)}" description]                                                                                                                                                                                                                                                                                                                                                                                                                                                                                                                                                         &  &                                &                                                                                                                                    & |[alias=p1p2]| \mathbf{0}                                                                                                      \\
                                                                                                                                                                        & |[alias=s2]|\Jasso(S_2,0) \arrow[ru, "{\tiny (P_2,0)}"' description] \arrow[ld, "{\tiny (0,P_1)}" description] &                         &  &                                                                                                                                                                                                                                                                                                                                                                                                                                                                                                                                                                                                                                                                                       &  &                                &                                                                                                                                    &                                                                                                                               \\
|[alias=s2Sp1]| \mathbf{0}                                                                                                                                          &                                                                                                                                    &                         &  &                                                                                                                                                                                                                                                                                                                                                                                                                                                                                                                                                                                                                                                                                       &  &                                &                                                                                                                                    &                                                                                                                               \\
                                                                                                                                                                        &                                                                                                                                    &                         &  &                                                                                                                                                                                                                                                                                                                                                                                                                                                                                                                                                                                                                                                                                       &  &                                &                                                                                                                                    &                                                                                                                               \\
|[alias=Sp1]|\Jasso(0,P_1) \arrow[uu, "{\tiny (S_2,0)}" description] \arrow[dddd, "{\tiny (0,S_2)}" description] &                                                                                                                                    &                         &  & |[alias=o]| \fgMod(\Lambda) \arrow[rrrr, "{\tiny{(P_1,0)}}" description] \arrow[rrrddd, "{\tiny {(S_1,0)}}" description] \arrow[dddd, "{\tiny {(0, P_2)}}" description] \arrow[llll, "{\tiny{(0,P_1)}}" description] \arrow[uuuu, "{\tiny {(P_2,0)}}" description] \arrow[llluuu, "{\tiny {(S_2,0)}}" description] \arrow[lluuuu, "{\tiny {(P_2\oplus S_2,0)}}" description,color=gray,opacity=0.4,near end] \arrow[rrrruuuu, "{\tiny{(P_1\oplus P_2,0)}}" description,color=gray,opacity=0.4,near end] \arrow[rrrrdd, "{\tiny{(P_1\oplus S_1,0)}}" description,color=gray,opacity=0.4,near end] \arrow[rrdddd, "{\tiny {(S_1, P_2)}}" description,color=gray,opacity=0.4,near end] \arrow[lllldddd, "{\tiny {(P_1\oplus P_2,0)}}" description,color=gray,opacity=0.4,near end] \arrow[lllluu, "{\tiny {(S_2, P_1)}}" description,color=gray,opacity=0.4,near end] &  &                                &                                                                                                                                    & |[alias=p1]| \Jasso(P_1,0) \arrow[uuuu, "{\tiny (S_2,0)}" description] \arrow[dd, "{\tiny (0,S_2)}" description] \\
                                                                                                                                                                        &                                                                                                                                    &                         &  &                                                                                                                                                                                                                                                                                                                                                                                                                                                                                                                                                                                                                                                                                       &  &                                &                                                                                                                                    &                                                                                                                               \\
                                                                                                                                                                        &                                                                                                                                    &                         &  &                                                                                                                                                                                                                                                                                                                                                                                                                                                                                                                                                                                                                                                                                       &  &                                &                                                                                                                                    & |[alias=p1s1]| \mathbf{0}                                                                                                   \\
                                                                                                                                                                        &                                                                                                                                    &                         &  &                                                                                                                                                                                                                                                                                                                                                                                                                                                                                                                                                                                                                                                                                       &  &                                & |[alias=s1]|\Jasso(S_1,0) \arrow[ld, "{\tiny (0,P_2)}"' description] \arrow[ru, "{\tiny (P_1,0)}" description] &                                                                                                                               \\
|[alias=Sp1Sp2]|\mathbf{0}                                                                                                                                 &                                                                                                                                    &                         &  & |[alias=Sp2]|\Jasso(0,P_2) \arrow[llll, "{\tiny (0,S_1)}"' description] \arrow[rr, "{\tiny (S_1,0)}" description]                                                                                                                                                                                                                                                                                                                                                                                                                                                                                                                &  & |[alias=s1Sp2]| \mathbf{0} &                                                                                                                                    &                                                                                                                              
\end{tikzcd}
\end{equation*}
Again, the parts that have the same colour/shading are to be identified. The faint arrows are morphisms given by $\tau$-rigid pairs with two indecomposable direct summands, and the solid arrows are given by indecomposable $\tau$-rigid pairs.
The picture space of $\Lambda$ is homeomorphic to a twice punctured torus.

\section{Ringel--Hall algebras}\label{sec:RingelHall}
\setcounter{subsection}{-1}

\subsection{}\label{Hall_intro}
We recall the construction of Ringel--Hall algebras of stack functions. We follow Joyce \cite{Joy07} and Bridgeland  \cite{Bri17}, who build on a long history of Hall algebras \cite{Ste01,Hal59,Rin90,Lus91,Nak94,KV00,FMV01} (we refer to Schiffmann \cite{Sch06} for a thorough introduction to the classical theory). In this section and henceforth, we assume the field $\field$ to be algebraically closed.

\subsection{} 
A \textit{pre-stack} over $\field$ is a contravariant 2-functor from the category of schemes over $\field$ (regarded as a 2-category by declaring the 2-morphisms to be identities) to the 2-category of groupoids. 
A \textit{stack} is a pre-stack satisfying some conditions we do not recite here (see \cite[Definition 2.10]{Gom01}  \cite[§2.3]{Joy06}). 
The 2-category of pre-stacks is simply a 2-functor category, and there is a full 2-subcategory of stacks over $\field$.

Any scheme over $\field$ can be regarded as a stack over $\field$, by composing its functor of points with the inclusion of the category of sets into the 2-category of groupoids. We will in particular regard $\mathrm{Spec}(\field)$ as a stack. A \textit{$\field$-point} in a stack $\mathscr{X}$ over $\field$ is a 2-isomorphism class of 1-morphisms $\mathrm{Spec}(\field)\to \mathscr{X}$. We denote the set of \textit{$\field$-points} in $\mathscr{X}$ by $\mathscr{X}(\field)$. Let $\mathrm{Iso}_{\field}(x)$ denote the \textit{geometric stabilizer group} of the $\field$-point $x$ in $\mathscr{X}$, namely the group of 2-endomorphisms of $x$.
A a stack $\mathscr{X}$ over $\field$ is an \textit{algebraic stack} if the following two conditions are met \cite[§3]{LMB99}:
\begin{enumerate}
	\item The diagonal map $\Delta_{\mathscr{X}}\colon \mathscr{X} \to \mathscr{X}\times \mathscr{X}$ is representable, quasi-compact, and separated.
	\item There exists a scheme $X$ and a smooth surjective morphism $X \to \mathscr{X}$.
\end{enumerate}
When $\mathscr{X}$ is algebraic, all geometric stabilizer groups are algebraic groups over $\field$. If all geometric stabilizer groups are affine, we say that $\mathscr{X}$ \textit{has affine geometric stabilizers}. 

\subsection{}\label{modstack}
We have fixed a finite-dimensional algebra $\Lambda$ over an algebraically closed field $\field$.
By a classical result of Gabriel, one can find a finite quiver $Q$ and an ideal $I$ of the path algebra $\field Q$, such that $\Lambda$ is Morita equivalent to $\field Q/ I$. For a scheme $X$ over $\field$, consider the groupoid in which the objects are pairs $(\mathcal{F},f)$, where $\mathcal{F}$ is a locally free $\mathcal{O}_X$-module of finite type and $f\colon \Lambda \to \mathrm{End}_{X}(\mathcal{F})$ is a homomorphism of $\field$-algebras, and a morphism from $(\mathcal{F}_1,f_1)$ to $(\mathcal{F}_2,f_2)$ is given by an isomorphism $\phi\colon\mathcal{F}_1\to \mathcal{F}_2$ of $\mathcal{O}_X$-modules such that the following square commutes for all $\lambda\in\Lambda$:
\begin{equation*}
	\begin{tikzcd}
		 \mathcal{F}_1 \arrow[r,"f_1(\lambda)"] \arrow[d,"\phi"] &  \mathcal{F}_1 \arrow[d,"\phi"] \\
		 \mathcal{F}_2 \arrow[r,"f_2(\lambda)"] &  \mathcal{F}_2
	\end{tikzcd}
\end{equation*}
One can define a prestack $\modstack$ by sending $X$ to this groupoid, and inducing morphisms using inverse images. This is an algebraic stack of locally finite type \cite[Theorem 8.2]{Joy06}  with affine geometric stabilizers (since the stabilizer group $\mathrm{Iso}_{\field}([X])$ is the group of invertible elements in the finite-dimensional $\field$-algebra $\mathrm{End}(X)$ \cite[§3.2]{JoyIII}). We will refer to it as the \textit{moduli stack} of $\Lambda$. 
The set of $\field$-points $\modstack(\field)$ is given by the set of isomorphism classes of finite-dimensional representations of the bound quiver $(Q,I)$, or equivalently by isomorphism classes of finite-dimensional $\Lambda$-modules.

\subsection{}\label{flags}
Let $m$ be a positive integer. We define the stack the $\modstack^{(m)}$ as the stack of diagrams of $\Lambda$-modules of the form
\begin{equation}\label{eq:Waldhausen}
	\begin{tikzcd}
		M_{[1,1]} \arrow[d,two heads]  \arrow[r,tail] \arrow[rd,phantom,"\square" description] & M_{[1,2]}  \arrow[r,tail] \arrow[d,two heads]\arrow[rd,phantom,"\square" description] & M_{[1,3]}  \arrow[r,tail] \arrow[d,two heads] &\cdots  \arrow[r,tail] & M_{[1,m]} \arrow[d,two heads] \\
		0 \arrow[r,tail] & M_{[2,2]}  \arrow[r,tail]  \arrow[d,two heads]\arrow[rd,phantom,"\square" description] & M_{[2,3]} \arrow[d,two heads]  \arrow[r,tail] & \cdots   \arrow[r,tail]& M_{[2,m]} \arrow[d,two heads] \\
		{} & 0\arrow[r,tail] & \ddots \arrow[d,two heads]  \arrow[r,tail]  \arrow[rd,phantom,"\square" description]& \ddots \arrow[r,tail]\arrow[d,two heads]  \arrow[rd,phantom,"\square" description] & \vdots \arrow[d,two heads] \\
		{} & {} &  0\arrow[r,tail] & M_{[m-1,m-1]}  \arrow[d,two heads] \arrow[r,tail] \arrow[rd,phantom,"\square" description] & M_{[m-1,m]} \arrow[d,two heads] \\
		{} & {} & {} & 0 \arrow[r,tail]& M_{[m,m]} 
	\end{tikzcd}
\end{equation}
where the horizontal $\Lambda$-homomorphisms are monomorphisms, the vertical morphisms are epimorphisms, and each commutative square is bicartesian. In particular, for every interval $[i,j]$ and $h\in[i,j-1]$, the sequence
\begin{equation*}
\begin{tikzcd}
	M_{[i,h]} \arrow[r,tail] & M_{[i,j]} \arrow[r,two heads] & M_{[h+1,j]} 
\end{tikzcd}
\end{equation*}
is short exact. The stacks $\modstack^{(m)}$ have affine geometric stabilizers \cite[Discussion following Definition 3.5]{Joy07}. Let $a_i^{} \colon \modstack^{(m)}\to \modstack$ be the restriction morphism sending the diagram in \eqref{eq:Waldhausen} to $M_{[i,i]}$ (on the diagonal), and let $b^{} \colon \modstack^{(m)}\to \modstack$ be the restriction morphism sending the diagram in \eqref{eq:Waldhausen} to $M_{[1,m]}$ (top right). 
 
 \subsection{}\label{SF}
Given be an algebraic stack $\mathscr{X}$ over $\field$ with affine geometric stabilizers, we will work with pairs $(\mathscr{R},\rho)$, where $\mathscr{R}$ is an algebraic stack over $\field$ of finite type with affine geometric stabilizers and $\rho$ is a 1-morphism from $\mathscr{R}$ to $\mathscr{X}$.  
Two such pairs $(\mathscr{R}_1,{\rho_1})$ and $(\mathscr{R}_2,{\rho_2})$ are \textit{equivalent} if there exists a 1-isomorphism $\mathscr{R}_1 \to^{\phi} \mathscr{R}_2$ and a 2-isomorphism between $\rho_2 \circ \phi$ and $\rho_1$. An equivalence class under this relation will be denoted $[(\mathscr{R},\rho)]$. In the free abelian group $F$ generated by these equivalence classes, consider the subgroup $K$ spanned by the expressions \ref{SFadd}--\ref{SFfib} below. 
\begin{enumerate}
	\item\label{SFadd} $[(\mathscr{R},\rho)] - [(\mathscr{Q},\rho\vert_\mathscr{Q})] - [(\mathscr{R}\setminus\mathscr{Q},\rho\vert_{\mathscr{R}\setminus\mathscr{Q}})]$,
where $\mathscr{Q}$ is an open substack  of $\mathscr{R}$. 
	\item\label{SFbij} $[(\mathscr{R}_1,\rho_1)] - [(\mathscr{R}_2,\rho_2)]$, if there exists a strictly commutative diagram of stacks
	\[
	\begin{tikzcd}
		\mathscr{R}_1 \arrow[rr,"\phi"] \arrow[rd,"\rho_1"']& & \mathscr{R}_2 \arrow[ld,"\rho_2"] \\
		 & \mathscr{X} &
	\end{tikzcd}
	\]
where $\phi$ induces an equivalence between the groupoids of $k$-points. 
	\item\label{SFfib} $[(\mathscr{Y}_1, \rho\circ \omega_1)] - [(\mathscr{Y}_2, \rho\circ \omega_2)]$, for every diagram 
	\[
	 \begin{tikzcd}[row sep=0em]
	 	\mathscr{Y}_1 \arrow[rd,"\omega_1", near start] && \\
		& \mathscr{R} \arrow[r,"\rho"] & \mathscr{X} \\
		\mathscr{Y}_2 \arrow[ru,"\omega_2"', near start] &&
	 \end{tikzcd}
	 \]
	where $\omega_1$ and $\omega_2$ are locally-trivial fibrations in the Zariski topology with isomorphic fibres.
\end{enumerate}
The \textit{Grothendieck group of stacks over $\mathcal{X}$} is then defined by $K(\mathbf{St}/\mathscr{X}) \defeq F/K$ (cf. \cite[Definition 3.6]{Bri11}). 

Let $K(\mathbf{St}/\field) \defeq K(\mathbf{St}/\mathrm{Spec}(\field))$. Using fibre products, one equips $K(\mathbf{St}/\field)$ with the structure of commutative ring. If $\mathscr{X}$ is as in the previous paragraph, one gives $K(\mathbf{St}/\mathscr{X})$ a $K(\mathbf{St}/\field)$-module structure by imposing 
\[  [ (\mathscr{Y}, \mathscr{Y}\to^{} \mathrm{Spec}(\field) )] \cdot [(\mathscr{R}, \mathscr{R}\to^{\rho} \mathscr{X})] \defeq  [(\mathscr{Y}\times \mathscr{R}, \mathscr{Y}\times \mathscr{R}\to^{\rho\circ \pi} \mathscr{X})], \]
where $\pi$ is the projection $\mathscr{Y}\times \mathscr{R}\to \mathscr{R}$.
The function
 \[\Upsilon \colon \{\text{isomorphism classes of quasi-projective $\field$-varieties}\} \to \Q(t) \]
assigning virtual Poincaré polynomials (we refer to Joyce \cite[Example 4.3]{Joy07motivic} for a definition), induces a ring homomorphism 
 \[\Upsilon \colon K(\mathbf{St}/\field) \to \Q(t). \]
In the following, $K_{\Upsilon}(\mathbf{St}/\mathscr{X})$ will denote the $\Q(t)$-vector space $\Q(t)\otimes_{K(\mathbf{St}/k)}K(\mathbf{St}/\mathscr{X})$.
 
 \subsection{}\label{HallAlg}
 We are now ready to recall the definition of the (completed) Ringel--Hall algebra of $\Lambda$. One defines an associative product on $K_{\Upsilon}(\mathbf{St}/\modstack)$ by setting \[[(\mathscr{X}_1,\rho_1)]\HallProd [(\mathscr{X}_2,\rho_2)] = [(\mathscr{Z},{b\circ \widehat{\rho}})], \]
 where $\widehat{\rho}$ appears in the following diagram, in which the square is homotopy cartesian \cite[§5.3]{Bri17} \cite[Discussion following Definition 2.15]{Joy07} 
 \begin{equation}\label{eq:torsioncomp}
	\begin{tikzcd}[]
		\mathscr{Z} \arrow[r,"{\widehat{\rho}}"] \arrow[d]& \modstack^{(2)} \arrow[r,"b"] \arrow[d,"{(a_1,a_2)}"] & \modstack \\
		\mathscr{X}_1 \times \mathscr{X}_2 \arrow[r,"\rho_1\times \rho_2"] & \modstack\times\modstack & {}
	\end{tikzcd}
	\end{equation}
The identity element in $K_{\Upsilon}(\mathbf{St}/\modstack)$ is given by $1_{\mathbf{0}} \defeq [(\modstack,\modstack_{\mathbf{0}}\hookrightarrow \modstack)]$, where $\modstack_0$ denotes the stack of zero $\Lambda$-modules.
We refer to $K_{\Upsilon}(\mathbf{St}/\modstack)$ as the \textit{Ringel--Hall algebra (of stack functions)} of $\Lambda$. It will be denoted $H(\Lambda)$ henceforth. 

Let $N^{\oplus}$ denote the free abelian group $\mathbb{N}^n$, where $n$ is the \textit{rank} of $\Lambda$ (i.e. the number of simple $\Lambda$-modules up to isomorphism). The $\Q(t)$-algebra $H(\Lambda)$ is then $N^{\oplus}$-graded:
\begin{equation*}
	H(\Lambda) = \bigoplus_{\mathbf{d}\in N^{\oplus}} H_{\mathbf{d}}(\Lambda),
\end{equation*}
where $H_{\mathbf{d}}(\Lambda) \defeq K_{\Upsilon}(\mathbf{St}/\modstack_{\mathbf{d}})$ and $\modstack_{\mathbf{d}}$ denotes the open substack of $\modstack$ parameterizing the $\Lambda$-modules with dimension vector $\mathbf{d}$. If $\delta$ denotes the summation map from $N^{\oplus}$ to $\mathbb{N}$, one turns $H(\Lambda)$ into a filtered $\Q(t)$-algebra
\begin{equation*}
	H(\Lambda) \supseteq H_{>0}(\Lambda) \supseteq H_{>1}(\Lambda) \supseteq \cdots, 
\end{equation*}
where $H_{> i}(\Lambda) \defeq \bigoplus_{\delta(\mathbf{d})>i} H_{\mathbf{d}}(\Lambda)$. The complete $\Q(t)$-algebra
\[ \widehat{H}(\Lambda) \defeq \varprojlim H(\Lambda)/H_{> i}(\Lambda), \]
is the \textit{completed Ringel--Hall algebra} of $\Lambda$. The subset $\widehat{G}\defeq 1_{\mathbf{0}} + \widehat{H}_{>0}(\Lambda)$ forms a group under $\HallProd$ (see \cite[§5.6]{Bri17}). It will play a major role in the following.

\subsection{}\label{Bri17.5.8}
Let $\mathscr{X} \to^{\rho} \modstack$ be a 1-morphism of stacks such that $\mathscr{X}$ has affine geometric stabilizers and $\mathscr{X}_{\mathbf{d}} \defeq \rho^{-1}(\modstack_{\mathbf{d}})$ is of finite type for all $\mathbf{d}\in N^{\oplus}$. 
Then the sequence 
\[ \left(\bigsqcup\limits_{\delta(\mathbf{d})\leq i} \mathscr{X}_{\mathbf{d}} \to^{\rho_{\leq i}} \modstack \right)_{i\geq 0} \]
determines an element of $\widehat{H}(\Lambda)$ (see \cite[§5.8]{Bri17}). It will be denoted $[(\mathscr{X},\rho)]$. It lies in $\widehat{G}$ precisely when $\mathscr{X}_{\mathbf{0}} \to^{\rho_{\mathbf{0}}} \modstack_{\mathbf{0}}$ is a 1-isomorphism of stacks.

\subsection{}\label{Bri12}
Let $\mathcal{X}$ and $\mathcal{Y}$ be full subcategories of $\fgMod(\Lambda)$. Assume that there exist open substacks $\modstack_{\mathcal{X}}$ and $\modstack_{\mathcal{Y}}$ of $\modstack$ whose sets of $\field$-points are the sets of isomorphism classes of objects $\mathcal{X}$ and $\mathcal{Y}$, respectively. By \Cref{Bri17.5.8}, these substacks determine group elements of $\widehat{G}$, which we denote $1_{\mathcal{X}}$ and $1_{\mathcal{Y}}$, respectively.
	By definition, the product $1_{\mathcal{X}} \HallProd 1_{\mathcal{Y}}$ is given by $[(\mathscr{Z},b\circ\widehat{\iota})]$, where the square in the diagram below is homotopy cartesian
	\begin{equation}\label{eq:torsioncomp}
	\begin{tikzcd}
		\mathscr{Z} \arrow[r,"\widehat{\iota}"] \arrow[d]& \modstack^{(2)} \arrow[r,"b"]\arrow[d,"{(a_1,a_2)}"] & \modstack \\
		\modstack_{\mathcal{X}} \times \modstack_{\mathcal{Y}} \arrow[r,hook] & \modstack\times\modstack & {}
	\end{tikzcd}
	\end{equation}
	Since the morphism $(a_1,a_2)$ satisfies the iso-fibration property of \cite[Lemma A.1]{Bri12}, the groupoid $\mathscr{Z}(\field)$ can be described in the following way (cf. \cite[Lemma 4.1]{Bri11}): 
	\begin{itemize}
		\item The objects are isomorphism classes short exact sequences of in $\fgMod(\Lambda)$ of the form 
	\begin{equation}\label{eq:torsionseq}
	\begin{tikzcd}
		 X  \arrow[r,tail] & E \arrow[r,two heads] & Y,
	\end{tikzcd}
	\end{equation}
	such that $X$ and $Y$ lie in the subcategories $\mathcal{X}$ and $\mathcal{Y}$, respectively.
	\item The morphisms are isomorphisms of short exact sequences.
	\end{itemize}
	  The function $(b\circ \widehat{\iota})(k)$ sends the isomorphism class of short exact sequences in \eqref{eq:torsionseq} to the isomorphism class of $E$. If $\mathcal{Y}\subseteq\mathcal{X}^{\perp}$, every object $E$ in $\mathcal{X}\VerdierProd\mathcal{Y}$ fits into a unique short exact sequence of the form displayed in \eqref{eq:torsionseq}. In this case, the 1-morphism $b\circ \widehat{\iota}$ induces an equivalence on the level of $k$-points. The relations \ref{SFbij} in \Cref{SF} identify $[(\mathscr{Z},b\circ \widehat{\iota})]$ with $[(\modstack,\mathrm{id}_{\modstack})]$.

\subsection{Lemma}\label{§2Lemma} 
Let $\mathcal{W}$ be a $\tau$-perpendicular subcategory of $\fgMod(\Lambda)$. 
\begin{enumerate}
	\item\label{§2Lemma1}  There is an open substack $\modstack_{\mathcal{W}}$ of $\modstack$ such that $\modstack_{\mathcal{W}}(\field)$ is the set of isomorphism classes of objects in $\mathcal{W}$. It determines an element $1_{\mathcal{W}}$ of $\widehat{G}$.
	\end{enumerate}
	Let $\mathcal{T}$ be a functorially finite torsion class of $\mathcal{W}$.
	\begin{enumerate}
	\setcounter{enumi}{1}
	\item\label{§2Lemma2} There is an open substack $\modstack_{\mathcal{T}}\subseteq \modstack$ such that $\modstack_{\mathcal{T}}(\field)$ is the set of isomorphism classes of objects in $\mathcal{T} \cap \mathcal{W}$. It determines a group element $1_{\mathcal{T}}$ of $\widehat{G}$. A similar statement holds for the functorially finite torsion-free class $\mathcal{T}^{\perp}$.
	\item\label{§2Lemma4} We have the \textit{torsion pair identity} $1_{\mathcal{W}} = 1_{\mathcal{T}} \HallProd 1_{\mathcal{T}^{\perp}}$ in $\widehat{G}$.
	\item\label{§2Lemma3.5} Let $\mathcal{U}$ be a functorially finite torsion class in $\mathcal{W}$ such that $\mathcal{U}\subseteq \mathcal{T}$. Then the group element $1_{[\mathcal{U},\mathcal{T}]}\defeq 1_{\mathcal{U}}^{-1}\HallProd 1_{\mathcal{T}}$ is determined by an open substack $\modstack_{[\mathcal{U},\mathcal{T}]}$ of $\modstack$ such that $\modstack_{[\mathcal{U},\mathcal{T}]}(\field)$ is the set of isomorphism classes of objects in $\mathcal{T}\cap\mathcal{U}^{\perp}\cap \mathcal{W}$.
	\item\label{§2Lemma3} Let $\mathcal{U}_1 \subseteq \mathcal{T}_1$ and $\mathcal{U}_2 \subseteq \mathcal{T}_2$ be two intervals in the poset of functorially finite torsion classes in $\mathcal{W}$. Then $1_{[\mathcal{U}_1,\mathcal{T}_1]}=1_{[\mathcal{U}_2,\mathcal{T}_2]}$ precisely when $\mathcal{T}_1\cap \mathcal{U}_1^{\perp}\cap \mathcal{W} = \mathcal{T}_2\cap \mathcal{U}_2^{\perp}\cap \mathcal{W}$. In particular, if $\mathcal{T}$ and $\mathcal{U}$ are distinct torsion classes in $\mathcal{W}$, then $1_{\mathcal{T}}$ and $1_{\mathcal{U}}$ are distinct elements of $\widehat{G}$.
\end{enumerate}
\begin{proof}
	The claims in \ref{§2Lemma1} and \ref{§2Lemma2} will be deduced from more a general result. We follow Bridgeland \cite[§6]{Bri17}. Observe that the proof is valid over any algebraically closed field $\field$, and not just the complex numbers.
	
	Let $M$ denote the dual abelian group of the Grothendieck group $K_0(\fgMod(\Lambda))$, and set $M_{\R}\defeq M\otimes_{\Z} \R$. Given $\theta\in M_{\R}$, one defines a stability function $Z_{\theta}: K_0(\fgMod(\Lambda)) \to \C$ by $Z_{\theta}([X])= -\theta([X]) + i \ell(X)$, where $\ell(X)$ denotes the length of the $\Lambda$-module $X$. Any $\Lambda$-module admits a unique Harder--Narasimhan filtration \cite[Theorem 3]{Rud97} with respect to $Z_\theta$. For an interval $I\subseteq [0,1]$, let $\mathcal{P}_{\theta}(I)$ denote the full subcategory of $\fgMod(\Lambda)$ spanned by the $\Lambda$-modules whose Harder--Narasimhan filtration with respect to $Z_\theta$ have phases in $I$. Bridgeland shows that there is an open substack $\modstack_{\theta}(I)\subseteq \modstack$ such that $\modstack_{\theta}(I)(\field)$ is the set of isomorphism classes of objects in $\mathcal{P}_{\theta}(I)$ \cite[§6.3]{Bri17}. Moreover, it determines an element $1_{\theta}(I)$ of $\widehat{G}$. 

 	We now prove \ref{§2Lemma1}. Since $\mathcal{W}$ is $\tau$-perpenducular, there exists a $\tau$-rigid pair $(M,P)$ such that $\mathcal{W}=\Jasso(M,P)$. Decomposing into indecomposables, we write $M\simeq \bigoplus_{i=1}^m M_i$ and $P\simeq \bigoplus_{j=m+1}^n P_j$. 
	Let $\theta_{(M,P)}$ denote the element of $M_{\R}$ defined by
	\begin{equation}\label{eq:thetaMP}
		\theta_{(M,P)}([X]) \defeq \sum_{i=1}^m \Hom_{\Lambda}(M_i,X) + \sum_{i=1}^m \Hom_{\Lambda}(X,\tau M_i) - \sum_{j=m+1}^n \Hom_{\Lambda}(P_j,X) \in \R. 
	\end{equation}
	Then $\theta_{(M,P)}$ lies in the interior of the $\bfg$-vector cone of $(M,P)$ \cite[Proof of Lemma 3.12]{BST19}. Consequently, the objects in $\mathcal{W}$ are exactly the $\theta_{(M,P)}$-semistable $\Lambda$-modules \cite[Theorem 3.14]{BST19}. Using the notation of the previous paragraph, this shows that $\mathcal{W}=\mathcal{P}_{\theta_{(M,P)}}(\{{1\over 2}\})$, so we may conclude our proof of \ref{§2Lemma1}.
	
	For the remainder of the proof, we may assume that $\mathcal{W}=\fgMod(\Lambda)$. Indeed, any open substack of $\modstack_{\mathcal{W}}$ is an open substack of $\modstack$, and $\mathcal{W}$ is equivalent to a module category.
	
	We move on to prove \ref{§2Lemma2}, setting $\mathcal{W}=\fgMod(\Lambda)$. If $\mathcal{T}$ is a functorially finite torsion class of $\fgMod(\Lambda)$, there exists a support $\tau$-tilting pair $(M,P)$ such that $\mathcal{T}=\fac(M)$ \cite[Theorem 2.7]{AIR14}. Let $\theta_{(M,P)}$ be as defined in \eqref{eq:thetaMP}. Since $\theta_{(M,P)}$ lies in the interior of the $\bfg$-vector cone of $(M,P)$, we have $\mathcal{T}=\mathcal{P}_{\theta_{(M,P)}}({1\over 2},1)$ \cite[Lemma 6.6]{Bri17} \cite[Remark 3.28]{BST19}, which proves our claim for functorially finite torsion classes. By showing that $\mathcal{T}^{\perp}=\mathcal{P}_{\theta_{(M,P)}}(0,{1\over 2}]$, one proceeds dually to prove the dual claim for functorially finite torsion-free classes.
	
	We now prove \ref{§2Lemma4}, setting $\mathcal{W}=\fgMod(\Lambda)$. In our proof of \ref{§2Lemma2}, we showed that $\mathcal{T}= \mathcal{P}_{\theta_{(M,P)}}({1\over 2},1)$ and that $\mathcal{T}= \mathcal{P}_{\theta_{(M,P)}}(0,{1\over 2}]$. Since the interval $(0,1)$ is a disjoint union $(0,1) = (0,{1\over 2}] \sqcup ({1\over 2},1)$, Bridgeland provides a decomposition \cite[Proposition 6.4]{Bri17} 
	\[1_{\fgMod(\Lambda)} = \mathcal{P}_{\theta_{(M,P)}}(0,1) = \mathcal{P}_{\theta_{(M,P)}}({1\over 2},1) \HallProd \mathcal{P}_{\theta_{(M,P)}}(0,{1\over 2}] = 1_{\mathcal{T}}\HallProd 1_{\mathcal{T}^{\perp}},\]
	as claimed.
	
	Next, we prove \ref{§2Lemma3.5}, setting $\mathcal{W}=\fgMod(\Lambda)$. The intersection $\modstack_{[\mathcal{U},\mathcal{T}]} \defeq \modstack_{\mathcal{T}} \times_{\modstack} \modstack_{\mathcal{U}^{\perp}}$ is an open substack of $\modstack$ whose $\field$-points is given by the set of isomorphism classes of objects in $\mathcal{T}\cap\mathcal{U}^{\perp}$. Using \Cref{Bri12}, one shows that the group element $g$ determined by $\modstack_{[\mathcal{U},\mathcal{T}]}$ satisfies $1_{\mathcal{T}} = 1_{\mathcal{U}} \HallProd g$. This is enough to prove our claim.
		
	Finally, we prove \ref{§2Lemma3}, again setting $\mathcal{W}=\fgMod(\Lambda)$. If we have $\mathcal{T}_1\cap \mathcal{U}_1^{\perp} = \mathcal{T}_2\cap \mathcal{U}_2^{\perp}$, then the open substacks determining these subcategories coincide, whence the equality $1_{[\mathcal{U}_1,\mathcal{T}_1]}=1_{[\mathcal{U}_2,\mathcal{T}_2]}$ follows. Conversely, suppose that  $1_{[\mathcal{U}_1,\mathcal{T}_1]}=1_{[\mathcal{U}_2,\mathcal{T}_2]}$. By symmetry, it suffices to show that $\mathcal{T}_1\cap \mathcal{U}_1^{\perp} \subseteq \mathcal{T}_2\cap \mathcal{U}_2^{\perp}$. The assumption implies that $1_{\mathcal{T}_1}=1_{\mathcal{U}_1}\HallProd 1_{[\mathcal{U}_2,\mathcal{T}_2]}$. It follows from \Cref{Bri12} that for every object $T_1\in\mathcal{T}_1$, there exists a short exact sequence $\begin{tikzcd} U_1 \arrow[r,tail] & T_1 \arrow[r,two heads] & F_1, \end{tikzcd}$ with $U_1\in \mathcal{U}_1$ and $F_1\in \mathcal{T}_2\cap \mathcal{U}_2^{\perp}$. In particular, if $T_1\in  \mathcal{T}_1\cap \mathcal{U}_1^{\perp}$, we deduce that $T_1\in  \mathcal{T}_2\cap \mathcal{U}_2^{\perp}$, which was what we set out to show.
	 \end{proof}

\section{Faithful group functors}
\setcounter{subsection}{-1}

\subsection{} In this section, we use Ringel--Hall algebras to construct a faithful functor from the $\tau$-cluster morphism category $\tcmc{\Lambda}$ into the group $\widehat{G}$ (as defined in \Cref{HallAlg}). This will simplify Igusa's sufficient conditions for when the picture space $\picspace{\Lambda}$ is locally $\mathrm{CAT}(0)$, at least under our assumption that $\field$ is an algebraically closed field (and more generally when $k$ is perfect \cite[Corollary 4.13]{BHK25}). The finite-dimensional $\field$-algebra $\Lambda$ remains fixed. 

\subsection{Theorem}\label{thm:2main}
Let $\tcmc{\Lambda}$ denote the $\tau$-cluster morphism category of $\Lambda$ and let $\picgroup{\Lambda}$ denote the picture group of $\Lambda$.
\begin{enumerate}
\item\label{thm:2main1} Regard the group $\widehat{G}$ as a groupoid with a single object. We may define a faithful functor $\eta\colon \tcmc{\Lambda} \to \widehat{G}$ by sending a morphism $\mathcal{W}\to^{(M,P)}\mathcal{V}$ to $1_{\fac_{\mathcal{W}}(M)}$ in $\widehat{G}$.
\item\label{thm:2main2} The generators $g_B$ of the picture group $\picgroup{\Lambda}$, where $B$ runs over the f-bricks in $\fgMod(\Lambda)$, are pair-wise distinct.
\end{enumerate}
\begin{proof}
We first prove \ref{thm:2main1}. 
The fact that $\eta$ gives a well-defined assignment of objects is trivial, and the well-definedness of the assignment of morphisms follows from \nref{§2Lemma}\ref{§2Lemma2}.  An identity morphism in $\tcmc{\Lambda}$ is of the form $\mathcal{W}\to^{(0,0)}\mathcal{W}$, whence it is clear from the construction that $\eta$ sends identity morphisms in $\tcmc{\Lambda}$ to the identity morphism $1_{\mathbf{0}}$ in $\widehat{G}$. Consider a pair of composable morphisms in $\tcmc{\Lambda}$, as displayed along the top row below, as well as their composite displayed by the bent arrow. 
\begin{equation*}
	\begin{tikzcd}[column sep=5em]
\mathcal{W} \arrow[r, "{(M,P)}"] \arrow[rr, "{E^{-1}_{(M,P)}(N,Q)}"', bend right=20] & \mathcal{V} \arrow[r, "{(N,Q)}"] & \mathcal{V}',
\end{tikzcd}
\end{equation*}
To show that $\eta$ is a functor, it remains to show that the composite ${E^{-1}_{(M,P)}(N,Q)}$ is sent to the product $1_{\fac_{\mathcal{W}}(M)}\HallProd 1_{\fac_{\mathcal{V}}(N)}$ in $\widehat{G}$. Since this is a direct consequence of \nref{lem:facMM'} and the fact that ${\fac_{\mathcal{V}}(N)} \subseteq {\fac_{\mathcal{W}}(M)}^{\perp}$, we move on to show that $\eta$ is faithful. Let $\calW_1$ and $\calW_2$ be objects in $\tcmc{\Lambda}$, and consider two distinct morphisms $(M_1,P_1), (M_2,P_2) \in \tcmc{\Lambda}(\calW_1,\calW_2)$. The group elements $1_{\fac_{\mathcal{W}}(M_1)}$ and $1_{\fac_{\mathcal{W}}(M_2)}$ are distinct elements of $\widehat{G}$ as a result of \nref{§2Lemma}\ref{§2Lemma3}, which shows that $\eta$ is faithful, as desired.

Finally, we prove \ref{thm:2main2}.
Let $B_1$ and $B_2$ be non-isomorphic f-bricks in $\fgMod(\Lambda)$, and let $g_1$ and $g_2$ denote the generators of the picture group $\picgroup{\Lambda}$ given by $B_1$ and $B_2$, respectively. For $i\in \{1,2\}$, the functor $\eta$ sends the morphism $\begin{tikzcd} \filt(B_i) \arrow[r,"\filt(B_i)"]& \mathbf{0} \end{tikzcd}$ to $1_{\filt(B_i)}$ and the morphism $\begin{tikzcd} \filt(B_i) \arrow[r,"\mathbf{0}"]& \mathbf{0} \end{tikzcd}$ to $1_{\mathbf{0}}$. Thus, the induced group homomorphism $\eta\colon \picgroup{\Lambda} \to \widehat{G}$ sends $g_i$ to $1_{\filt(B_i)}$. Since $\filt(B_i)$ is a $\tau$-perpendicular subcategory of $\fgMod(\Lambda)$, there exist functorially finite torsion classes $\mathcal{T}_i$ and $\mathcal{U}_i$ of $\fgMod(\Lambda)$ such that $\mathcal{U}_i \subseteq \mathcal{T}_i$ and $\mathcal{T}_i\cap\mathcal{U}_i^{\perp}=\filt(B_i)$ \cite[Theorem 1.1]{BH21}. The functor $\eta$ then sends $g_i$ to $1_{[\mathcal{U}_i,\mathcal{T}_i]}\in \widehat{G}$. Our claim now follows directly from \nref{§2Lemma}\ref{§2Lemma3}.
\end{proof}

\subsection{}\label{Igu14.3.4} 
The $\tau$-cluster morphism category $\tcmc{\Lambda}$ of $\Lambda$ is a cubical category, in the sense of Igusa \cite[§3.1]{Igu14} \cite[Theorem 2.14]{HI21}\footnote{The cited reference assumes that $\Lambda$ is {brick-finite}, but this assumption does not seem necessary for Hanson--Igusa in hindsight. Buan--Hanson later defined picture categories for arbitrary finite-dimensional $\field$-algebras \cite{BH21}.} (see also \cite[§6]{Bor24} and \cite[Theorem 3.16]{Kai23}). 
In \nref{cor:complastK_alg}, we will see that \nref{thm:2main} sometimes can be used to determine the homotopy type of picture spaces.
In greater generality, Igusa has shown that a cubical category $\mathcal{C}$ is locally $\mathrm{CAT}(0)$ when the conditions \ref{I1}--\ref{I3} below are met \cite[Proposition 3.5]{Igu14}. These conditions thus ensure that the classifying space of $\mathcal{C}$ is an Eilenberg--MacLane space of type $\mathrm{K}(\pi,1)$.

We do not recall Igusa's definition of a cubical category here. However, we will mention that every morphism in a cubical category is assigned a non-negative integer, called the \textit{rank}. If the rank of $f$ is $r$ and the rank of $g$ is $s$, then the rank of $g\circ f$ is $r+s$ whenever this composite is defined. A \textit{first factor} of a morphism $f$ is a morphism $f_1$ of rank 1 such that $f=g\circ f_1$ for some morphism $g$. The notion of \textit{last factor} is defined dually. The axioms of cubical categories guarantee that every morphism admits a first factor and a last factor. 

\begin{enumerate}[label=(\textbf{I\arabic*})]
	\item\label{I1} Let $X$ be an object in $\mathcal{C}$ and let $\{X \to Y_j \}$ be a set of $\ell\geq 2$ morphisms in $\mathcal{C}$, all of rank $1$. This set forms the set of first factors of a rank $\ell$ morphism in $\mathcal{C}$ if and only if each subset of cardinality 2 forms the set of first factors of a rank 2 morphism in $\mathcal{C}$.
	\item\label{I2} 
	Let $X$ be an object in $\mathcal{C}$ and let $\{Y_j \to X \}$ be a set of $\ell\geq 2$ morphisms in $\mathcal{C}$, all of rank $1$.  This set forms the set of last factors of a rank $\ell$ morphism in $\mathcal{C}$ if and only if each subset of cardinality 2 forms the set of last factors of a rank 2 morphism in $\mathcal{C}$.
	\item\label{I3} There exists a faithful functor $\mathcal{C}\to\Gamma$, for some group $\Gamma$ (regarded as a groupoid).
\end{enumerate}

\subsection{Corollary}\label{cor:complastK_alg}
 (cf. \cite[Theorem 4.13]{HI21})
 Let $\Lambda$ be a finite-dimensional $\field$-algebra and regard the $\tau$-cluster morphism category $\tcmc{\Lambda}$ as a cubical category. The conditions \ref{I1} and \ref{I3} in \Cref{Igu14.3.4} hold. 
 Thus, the condition \ref{I2} guarantees that the picture space $\picspace{\Lambda}$ is locally $\mathrm{CAT}(0)$. In particular, the condition \ref{I2} is sufficient for $\picspace{\Lambda}$ to be an Eilenberg--MacLane space of type $\mathrm{K}(\pi,1)$.
 \begin{proof}
 	The condtion \ref{I3} immediately follows from \nref{thm:2main}. The condition \ref{I1} is straightforward to verify for $\tcmc{\Lambda}$. Indeed, it is equivalent to the fact that a direct sum of indecomposable $\tau$-rigid pairs $(M_i,P_i)$ is $\tau$-rigid if and only if every direct sum of two of them is $\tau$-rigid. 
 \end{proof}
 
 \subsection{}
 Classes of examples where \ref{I2} holds include representation-finite hereditary $\field$-algebras \cite[Theorem 3.1]{IT17}, Nakayama $\field$-algebras \cite[Theorem 4.16]{HI21}, Nakayama-like $\field$-algebras \cite[Theorem 5.13(a)]{HI21p}, gentle $\field$-algebras in which every vertex has degree at most 2 \cite[Theorem 5.13(b)]{HI21p}, and all algebras of rank 3 or smaller \cite[Theorem 4.0.7]{BH22}. In these cases, one uses \nref{cor:complastK_alg} below to assert that the picture space $\picspace{\Lambda}$ is locally $\mathrm{CAT}(0)$, and thus a $\mathrm{K}(\pi,1)$ space.
Although \ref{I2} does not hold in general (a cluster tilted $\field$-algebra of type $A_4$ \cite[Future work]{HI21p} and a preprojective algebra of type $A_n$ for $n\geq 4$ \cite[Theorem 7.1.1]{BH22} being counter-examples), there seems to be no known example of a finite-dimensional $\field$-algebra whose picture space is decidedly not a $\mathrm{K}(\pi,1)$. Moreover, all known examples where the picture space has been shown to be $\mathrm{K}(\pi,1)$ are examples where the picture space is locally $\mathrm{CAT}(0)$.

 \section{Green sequences as positive expressions in the picture group}\label{sec:green}
 \setcounter{subsection}{-1}
 
 \subsection{} Maximal green sequences originate in Keller's work on quantum dilogarithm identities \cite{Kel11green}. We refer to Keller--Demonet \cite{Kel20} for a survey on the topic. When generalizing from hereditary algebras to arbitrary algebras, one may encode (maximal) green sequences as chains of torsion classes \cite{BST19} (see also \cite{Nag13,Qiu15}). In order to apply results from previous sections, the field $\field$ remains algebraically closed, and $\Lambda$ is a finite-dimensional $\field$-algebra.

\subsection{}\label{green_prep}
Let $\ftors(\Lambda)$ denote the poset of functorially finite torsion classes in $\fgMod(\Lambda)$. By \textit{cover relation}, we will mean a relation $\mathcal{U} \subseteq \mathcal{T}$ in $\ftors(\Lambda)$
such that $\mathcal{U}\neq \mathcal{T}$ and there exists no functorially finite torsion class $\mathcal{T}'$ with $\mathcal{U} \subsetneq \mathcal{T}'\subsetneq \mathcal{T}$. We denote a cover relation by $\mathcal{U} \lessdot \mathcal{T}$. 
For each such cover relation, the full subcategory $\mathcal{T}  \cap \mathcal{U}^{\perp}\subseteq \fgMod(\Lambda)$ is equal to $\filt(B)$ for some brick $B$ \cite[§§1.2 and 2.2]{Asa20} \cite[Theorem 3.3(b)]{DIRRT17}. We refer to $B$ as the \textit{brick label} of the cover relation $\mathcal{U} \lessdot \mathcal{T}$.

\subsection{Lemma}\label{noHom1} 
A brick $B$ in $\fgMod(\Lambda)$ is an f-brick if and only if there exists a cover relation $\mathcal{U}\lessdot \mathcal{T}$ in $\ftors(\Lambda)$ having $B$ as its brick label. 
\begin{proof}
	  Buan--Hanson show that a wide subcategory $\mathcal{W}\subseteq \fgMod(\Lambda)$ is $\tau$-perpendicular if and only if there exist functorially finite torsion classes $\mathcal{G}_1$ and $\mathcal{G}_2$ such that $\mathcal{G}_2 \subseteq \mathcal{G}_1$ and $\mathcal{W}=\mathcal{G}_1\cap\mathcal{G}_2^{\perp}$ \cite[Theorem 1.1]{BH21}. Sufficiency now follows from the discussion in \Cref{green_prep}. To prove necessity, suppose that $B$ is an f-brick, i.e. that $\filt(B)$ is a $\tau$-perpendicular subcategory of $\fgMod(\Lambda)$. Using Buan--Hanson's result again, we write $\filt(B)=\mathcal{G}_1\cap\mathcal{G}_2^{\perp}$. Since the interval $[\mathcal{G}_2,\mathcal{G}_1]$ is isomorphic to the poset $\ftors(\filt(B))$ \cite[Theorem 1.4(1)]{AP22}, where we regard $\filt(B)$ a module category of a local $\field$-algebra, we deduce that $\mathcal{G}_2\subseteq \mathcal{G}_1$ must be a cover relation from the fact that $\ftors(\Gamma)$ contains exactly two elements whenever $\Gamma$ is a local $\field$-algebra.
\end{proof}

\subsection{Definition} A \textit{green sequence} in $\fgMod(\Lambda)$ is a finite chain of cover relations in $\ftors(\Lambda)$\footnote{Br\"{u}stle--Smith--Treffinger allow the torsion classes in this chain to be arbitrary \cite[Definition 4.8]{BST19}, but show that they are necessarily functorially finite \cite[Proposition 4.9]{BST19}.} of the form
\begin{equation}\label{eq:greenseq}
	\mathbf{0} =  \mathcal{T}_0 \lessdot \mathcal{T}_{1} \lessdot \cdots \lessdot \mathcal{T}_{\ell-1} \lessdot \mathcal{T}_{\ell} =\mathcal{T}.
\end{equation}
We say that $\mathcal{T}$ \textit{admits a green sequence} if a green sequence of the form displayed in \eqref{eq:greenseq} exists.
A green sequence is a \textit{maximal green sequence} if $\mathcal{T} = \fgMod(\Lambda)$. The $\field$-algebra $\Lambda$ is said to \textit{admit a maximal green sequence} if the torsion class $ \fgMod(\Lambda)$ admits a green sequence.

 \subsection{}\label{preLemmaGreen}
 Let $\mathcal{T}$ be a functorially finite torsion class in $\fgMod(\Lambda)$ admitting a green sequence of the form as shown in \eqref{eq:greenseq}. If the brick labelling of the cover relations reads $(B_1,\dots, B_{\ell-1},B_{\ell})$ from left to right, let $g_i$ denote the generator of the picture group $\picgroup{\Lambda}$ provided by the brick $B_i$ (see \eqref{eq:picgroup_pres_genloop}). Let $g_{\mathcal{T}} \defeq g_{1} \cdots g_{\ell-1}  g_{\ell}\in\picgroup{\Lambda}$.  
 
\subsection{Lemma}\label{lem:Green}
We keep the setup and notation of \Cref{preLemmaGreen}.
\begin{enumerate}
	\item\label{lem:Green1} Let $\eta\colon \picgroup{\Lambda} \to \widehat{G}$ be the group homomorphism induced by the functor $\eta\colon \tcmc{\Lambda} \to \widehat{G}$ defined in \nref{thm:2main}. For an arbitrary generator $g_i$ of $\picgroup{\Lambda}$, we have $\eta(g_i)=1_{\filt(B_i)}$. Furthermore, we have $\eta(g_{\mathcal{T}})= 1_{\mathcal{T}}$.
	\item\label{lem:Green2} The group element $g_{\mathcal{T}}\in\picgroup{\Lambda}$ is independent of the choice of green sequence.
\end{enumerate}
\begin{proof}
We first prove \ref{lem:Green1}. For each $B_i$, the cover relation $\mathcal{T}_{i-1} \lessdot \mathcal{T}_i$ is labelled by $B_i$. It was shown in the proof of \nref{thm:2main}\ref{thm:2main2} that $\eta(g_i)=1_{\filt(B_i)}=1_{[\mathcal{T}_{i-1},\mathcal{T}_i]}$. By the definition of the elements $1_{[\mathcal{T}_{i-1},\mathcal{T}_i]}$ given in \nref{§2Lemma}\ref{§2Lemma3.5}, we find that 
\begin{equation*}
\eta(g_{\mathcal{T}})  = \eta(g_1\cdots g_{\ell-1}g_{\ell}) = \eta(g_1)\HallProd \cdots \HallProd \eta(g_{\ell-1})\HallProd \eta(g_{\ell}) = 1_{[\mathbf{0},\mathcal{T}_{1}]}\HallProd \cdots \HallProd 1_{[\mathcal{T}_{\ell-2},\mathcal{T}_{\ell-1}]} \HallProd 1_{[\mathcal{T}_{\ell-1},\mathcal{T}_{\ell}]} = 1_{\mathcal{T}},
\end{equation*}
as claimed. 

The assertion in \ref{lem:Green2} is known to hold when $\Lambda$ is brick-finite (see \cite[Proposition 4.4]{HI21} \cite[Lemma 9-6.3]{Rea16}). We argue that it holds in general. Suppose that we are given the green sequence displayed in \eqref{eq:greenseq} as well as another green sequence of the form
\begin{equation*}
	\mathbf{0} =  \mathcal{T}_0 \lessdot \mathcal{T}'_{1} \lessdot \cdots \lessdot \mathcal{T}'_{\ell-1} \lessdot \mathcal{T}_{\ell} =\mathcal{T}.
\end{equation*}
We may assume that $\mathcal{T}_{1} \neq  \mathcal{T}'_{1}$ and that $\mathcal{T}_{\ell-1} \neq  \mathcal{T}'_{\ell-1}$. The two green sequences thus form a polygon in the poset $\ftors(\Lambda)$. From here, the reasoning of Hanson--Igusa \cite[§4.2]{HI21} applies, whence we conclude.
\end{proof}
 
 \subsection{Lemma}\label{lem:fill}
 Let $\mathcal{T}$ be a functorially finite torsion class in $\fgMod(\Lambda)$ and let $B_1, \dots, B_{\ell-1}$, and $B_{\ell}$ be bricks such that $1_{\filt(B_1)}\HallProd \cdots \HallProd 1_{\filt(B_{\ell -1})} \HallProd 1_{\filt(B_{\ell})}=1_{\mathcal{T}}$ in $\widehat{G}$. Extending \Cref{Bri12} from binary products to $\ell$-ary products of elements in $\widehat{G}$, we may claim that a $\Lambda$-module $X$ is in $\mathcal{T}$ precisely when there exists a diagram up to isomorphism in $\fgMod(\Lambda)$
	\begin{equation}\label{eq:Xdia}
	\begin{tikzcd}
		B'_1 \arrow[d,two heads]  \arrow[r,tail] \arrow[rd,phantom,"\square" description] & X_{[1,2]}  \arrow[r,tail] \arrow[d,two heads]\arrow[rd,phantom,"\square" description] & X_{[1,3]}  \arrow[r,tail] \arrow[d,two heads] &\cdots  \arrow[r,tail] &X \arrow[d,two heads] \\
		0 \arrow[r,tail] & B'_2 \arrow[r,tail]  \arrow[d,two heads]\arrow[rd,phantom,"\square" description] & X_{[2,3]} \arrow[d,two heads]  \arrow[r,tail] & \cdots   \arrow[r,tail]& X_{[2,m]} \arrow[d,two heads] \\
		{} & 0\arrow[r,tail] & \ddots \arrow[d,two heads]  \arrow[r,tail]  \arrow[rd,phantom,"\square" description]& \ddots \arrow[r,tail]\arrow[d,two heads]  \arrow[rd,phantom,"\square" description] & \vdots \arrow[d,two heads] \\
		{} & {} &  0\arrow[r,tail] & B'_{\ell-1} \arrow[d,two heads] \arrow[r,tail] \arrow[rd,phantom,"\square" description] & X_{[m-1,m]} \arrow[d,two heads] \\
		{} & {} & {} & 0 \arrow[r,tail]& B'_{\ell}
	\end{tikzcd}
\end{equation}
where $B'_{i}\in \filt(B_i)$ and the sqaures are bicartesian. 
\begin{enumerate}
	\item\label{lem:fill0} The bricks $B_1, \dots, B_{\ell-1}$, and $B_{\ell}$ are pair-wise distinct up to isomorphism.
	\item\label{lem:fill1} We have $\Hom_{\Lambda}(B_i,B_j)=0$ whenever $1\leq i < j \leq \ell$.
	\item\label{lem:fill2} Given an object $X\in \mathcal{T}$, we have $X\in {^{\perp}B_{\ell}}$ precisely when $B'_{\ell}$ in \eqref{eq:Xdia} is zero.
\end{enumerate}
\begin{proof}
We first prove \ref{lem:fill0}. The element $1_{\filt(B_1)}\HallProd \cdots \HallProd 1_{\filt(B_{\ell -1})} \HallProd 1_{\filt(B_{\ell})}$ in $\widehat{G}$ is determined the 1-morphism $\modstack_{\mathcal{T}} \embed \modstack$, where $\modstack_{\mathcal{T}}$ is the open substack determined by $\mathcal{T}$. Consider the induced map $\modstack_{\mathcal{T}}(k)\embed \modstack(k)$. If there were repetitions among $B_1, \dots, B_{\ell-1}$, and $B_{\ell}$ up to isomorphism, this map would not be injective. Since it is, we have proved \ref{lem:fill0}.

We now prove \ref{lem:fill1}. Let $f$ be a $\Lambda$-homomorphism from $B_i$ to $B_j$ where $1\leq i<j\leq \ell$. It factors through its image: 
\begin{equation*}
	\begin{tikzcd}
		B_i \arrow[r,two heads] & \mathrm{Im}(f) \arrow[r,tail] & B_j.
	\end{tikzcd}
\end{equation*}
We assume that $ \mathrm{Im}(f) $ in non-zero in order to derive a contradiction. Since $\mathrm{Im}(f)$ is a factor of $B_i$, and $B_i \in \mathcal{T}$, we have that $\mathrm{Im}(f)\in \mathcal{T}$. We may thus form a diagram as in \eqref{eq:Xdia}, setting $X=\mathrm{Im}(f)$. Let $b$ denote the maximal $x$ such that $B'_x$ is non-zero. Since $B_b$ is the only simple object in $\filt(B_b)$ \cite[§1.2]{Rin76}, there is an epimorphism from $B'_\ell$ to $B_\ell$.
Since we have a composite epimorphism
\begin{equation*}
\begin{tikzcd}
	B_{i} \arrow[r,two heads] & \mathrm{Im}(f) \arrow[r,two heads] & B_b'  \arrow[r,two heads] &  B_b,
\end{tikzcd}
\end{equation*}
the brick $B_b$ is contained in the smallest torsion class containing $B_i$. Since this torsion class contains a unique brick up to isomorphism \cite[Proposition 3.1]{BTZ19}, we have shown that $B_b\simeq B_i$, whence $b=i$ by \ref{lem:fill0} above. A dual argument shows that the minimal $x$ such that $B'_x$ is non-zero is equal to $j$. The inequalities $i<j$ and $i\geq j$ cannot hold simultaneously. This contradiction derived from the assumption that $ \mathrm{Im}(f) $ is non-zero shows that $f=0$, and thus that $\Hom_{\Lambda}(B_i,B_j)=0$ since $f$ was arbitrarily chosen.

Finally, we prove \ref{lem:fill2}. If $B'_{\ell}$ is non-zero, we use that $B_\ell$ is the only simple object in the wide subcategory $\filt(B_\ell)$ \cite[§1.2]{Rin76} to find an epimorphism from $B'_\ell$ to $B_\ell$. In \eqref{eq:Xdia}, there is an epimorphism from $X$ to $B'_{\ell}$. Composing the two epimorphisms shows that $X\not\in {^{\perp}B_{\ell}}$. Conversely, if $B'_{\ell}$ is zero, then $X$ is an iterated extension of the $\Lambda$-modules $B'_i$ for $1\leq i\leq \ell -1$. Since $B'_i\in\filt(B_i)$, it follows that $X$ is an iterated extension of the bricks $B_i$ for $1\leq i\leq \ell -1$. By \ref{lem:fill1}, each $B_i$ is in ${^{\perp}B_{\ell}}$, whence we conclude that $X\in{^{\perp}B_{\ell}}$.
\end{proof}
 
  \subsection{Theorem}\label{thm:positive_expr}
 Let $\mathcal{T}$ be a functorially finite torsion class in $\fgMod(\Lambda)$. Given a green sequence
 \begin{equation}\label{eq:greenseqThm}
    	\mathbf{0} =  \mathcal{T}_0 \lessdot \mathcal{T}_{1} \lessdot \cdots \lessdot \mathcal{T}_{\ell-1} \lessdot \mathcal{T}_{\ell} = \mathcal{T}
\end{equation}
in $\fgMod(\Lambda)$, with brick labels reading $(B_{1}, \dots , B_{\ell-1},B_{\ell})$ from left to right, let $g_i$ denote the generator of the picture group $\picgroup{\Lambda}$ provided by the brick $B_i$ (see \eqref{eq:picgroup_pres_genloop}). 
Mapping the green sequence in \eqref{eq:greenseqThm} to the word $g_{1}\cdots g_{\ell-1} g_{\ell}$ of generators of $\picgroup{\Lambda}$ defines a bijection from the first to the second of the following sets:\footnote{Note that we may also define this bijection when $\mathcal{T}$ does not admit a green sequence. Both sets listed will then be empty.}
\begin{enumerate}
	\item green sequences in $\fgMod(\Lambda)$ in of the form displayed in \eqref{eq:greenseqThm},
	\item positive expressions in the picture group $\picgroup{\Lambda}$ (i.e. a word of generators) whose product is equal to $g_{\mathcal{T}}$.
\end{enumerate}
In particular, we have a bijection from the first to the second of the following sets:
\begin{enumerate}
	\item maximal green sequences in $\fgMod(\Lambda)$,
	\item positive expressions in the picture group $\picgroup{\Lambda}$ whose product is equal to $g_{\fgMod(\Lambda)}$.
\end{enumerate}
\begin{proof}
	The fact that our map is well-defined is a direct consequence of \nref{lem:Green}\ref{lem:Green2}. Injectivity is a direct result of the following fact: Given a functorially finite torsion class $\mathcal{T}$ and an f-brick $B$, there is at most one cover relation $\mathcal{U} \lessdot  \mathcal{T}$ labelled by $B$. This holds by a result of Asai--Pfeifer \cite[Theorem 6.7]{AP22}, asserting that the brick labelling determines a bijection from the set of cover relations under $\mathcal{T}$ to the set of simple objects the following wide subcategory of $\fgMod(\Lambda)$ (see \cite[Proposition 3.3]{MS17}) up to isomorphism:
\begin{equation*}
\WL(\mathcal{T})\defeq \{X \in \mathcal{T} \sth \forall(g\colon Y\to X)\in \mathcal{T} \colon \ker(g) \in \mathcal{T}\}.
\end{equation*}
		
	We conclude the proof by showing surjectivity. Let $\mathcal{T}$ be a functorially finite torsion class in $\fgMod(\Lambda)$ admitting a green sequence, and let $g_{1}g_{2}\dots g_{\ell-1} g_{\ell}$ be an arbitrary positive expression in $\picgroup{\Lambda}$ whose product is equal to $g_{\mathcal{T}}$, where the $g_i$ denotes the generator of the picture group $\picgroup{\Lambda}$ provided by the brick $B_i$.
	We show that there exists a green sequence
	\begin{equation*}
	\mathbf{0} =  \mathcal{T}'_0 \lessdot \mathcal{T}'_{1} \lessdot \cdots \lessdot \mathcal{T}'_{\ell-1} \lessdot \mathcal{T}'_{\ell} = \mathcal{T},
\end{equation*}
	with brick labels $(B_{1},\dots,B_{\ell-1}, B_{\ell})$ reading left to right. It suffices to find a cover relation $\mathcal{U} \lessdot \mathcal{T}$ labelled by $B_{\ell}$, since we may then set $\mathcal{T}_{\ell-1}=\mathcal{U}$ and iterate the argument. Consider the torsion class $\mathcal{U}=\mathcal{T}\cap {^{\perp}B_{\ell}}$ of $\fgMod(\Lambda)$. It does not contain the indecomposable $\Lambda$-module $B_{\ell}$. If we can prove that all proper factors of $B_{\ell}$ are in $\mathcal{U}$, and that $B_{\ell}$ is a factor of each object in $\mathcal{T}\setminus\mathcal{U}$, we are done \cite[Proposition 2.4]{BTZ19}. 
	
	By \nref{lem:Green}\ref{lem:Green1}, we have $1_{\filt(B_1)}\HallProd \cdots \HallProd 1_{\filt(B_{\ell -1})} \HallProd 1_{\filt(B_{\ell})}=1_{\mathcal{T}}$ in $\widehat{G}$. Hence, we may apply the assertions in \nref{lem:fill}. Suppose that $X$ is a proper factor of $B_{\ell}$. Since $B_{\ell}\in \mathcal{T}$ and $\mathcal{T}$ is a torsion class, we have $X\in \mathcal{T}$, so $X$ admits a diagram as shown in \eqref{eq:Xdia}. If $B'_\ell$ is non-zero, consider the composite $\Lambda$-homomorphism
	\begin{equation*}
		\begin{tikzcd}
			B_{\ell} \arrow[r,two heads] & X \arrow[r,two heads] & B'_{\ell} \arrow[r, two heads] & B_{\ell}, 
		\end{tikzcd}
	\end{equation*}
rendering $B_{\ell}$ a proper factor of itself. Since this is an impossibility, it follows from \nref{lem:fill}\ref{lem:fill2} that $X\in \mathcal{U}$.

Now, let $X \in \mathcal{T}\setminus\mathcal{U}$. By \nref{lem:fill}\ref{lem:fill2} again, we have that $B'_{\ell}$ in \eqref{eq:Xdia} is non-zero. The existence of the composite epimorphism
	\begin{equation*}
		\begin{tikzcd}
 X \arrow[r,two heads] & B'_{\ell} \arrow[r, two heads] & B_{\ell}
		\end{tikzcd}
	\end{equation*}
shows that $B_{\ell}$ indeed is a factor of $X$. This completes our proof. 
	\end{proof}

\bibliographystyle{ourIEEEstyle.bst}
\bibliography{0Aus}

\end{document}